\documentclass[letterpaper,12pt]{amsart}
\usepackage{amssymb}
\usepackage{xypic}
\xyoption{all}
\xyoption{arc}
\newtheorem{thm}{Th\'eor\`eme}[section]
\newtheorem{lemme}[thm]{Lemme}
\newtheorem{cor}[thm]{Corollaire}
\newtheorem{prop}[thm]{Proposition}
\newtheorem{defi}[thm]{D\'efinition}
\theoremstyle{definition}
\newtheorem{rem}[thm]{Remarque}

\renewcommand{\thequestion}{}
\def\lhd{\vartriangleleft}
\def\Gbar{\bar{G}}
\def\Cbar{\bar{C}}
\def\Hbar{\bar{H}}
\def\Mbar{\bar{M}}
\def\Nbar{\bar{N}}
\def\ebar{\bar{e}}
\def\fbar{\bar{f}}
\def\Gp{{\mathfrak{p}}}
\def\Gg{{\mathfrak{g}}}
\def\Gm{{\mathfrak{m}}}
\def\Gn{{\mathfrak{n}}}
\def\GI{{\mathfrak{I}}}
\def\GA{{\mathfrak{A}}}
\def\Sn{{\mathfrak{S}}}
\def\Db{{\mathcal D}^b}
\def\CC{{\mathcal{C}}}
\def\CS{{\mathcal{S}}}
\def\CP{{\mathcal{P}}}
\def\CO{{\mathcal{O}}}
\def\CE{{\mathcal{E}}}
\def\CL{{\mathcal{L}}}
\def\CN{{\mathcal{N}}}
\def\CM{{\mathcal{M}}}
\def\CH{{\mathcal{H}}}
\def\CI{{\mathcal{I}}}
\def\CJ{{\mathcal{J}}}
\def\CV{{\mathcal{V}}}
\def\CA{{\mathcal{A}}}
\def\CF{{\mathcal{F}}}
\def\CG{{\mathcal{G}}}
\def\CT{{\mathcal{T}}}
\def\CD{{\mathcal{D}}}
\def\CU{{\mathcal{U}}}
\def\CW{{\mathcal{W}}}
\def\BC{{\mathbf{C}}}
\def\BR{{\mathbf{R}}}
\def\BF{{\mathbf{F}}}
\def\BZ{{\mathbf{Z}}}
\def\BL{{\mathbf{L}}}
\def\BA{{\mathbf{A}}}
\def\eps{\varepsilon}
\def\Tr{\operatorname{Tr}\nolimits}
\def\ab{\operatorname{ab}\nolimits}
\def\Id{\operatorname{Id}\nolimits}
\def\Sc{\operatorname{Sc}\nolimits}
\def\Br{\operatorname{Br}\nolimits}
\def\GL{\operatorname{GL}\nolimits}
\def\SL{\operatorname{SL}\nolimits}
\def\Trd{\operatorname{Trd}\nolimits}
\def\Aut{\operatorname{Aut}\nolimits}
\def\Ass{\operatorname{Ass}\nolimits}
\def\Out{\operatorname{Out}\nolimits}
\def\Inn{\operatorname{Inn}\nolimits}
\def\Hom{\operatorname{Hom}\nolimits}
\def\Homb{\operatorname{\overline{Hom}}\nolimits}
\def\Ext{\operatorname{Ext}\nolimits}
\def\Tor{\operatorname{Tor}\nolimits}
\def\Mod{\operatorname{mod}\nolimits}
\def\proj{\operatorname{proj}\nolimits}
\def\perm{\operatorname{perm}\nolimits}
\def\sets{\operatorname{sets}\nolimits}
\def\Stmod{\operatorname{Stmod}\nolimits}
\def\good{\operatorname{good}\nolimits}
\def\im{\operatorname{im}\nolimits}
\def\krull{\operatorname{Krulldim}\nolimits}
\def\ram{\operatorname{ram}\nolimits}
\def\End{\operatorname{End}\nolimits}
\def\Ind{\operatorname{Ind}\nolimits}
\def\Res{\operatorname{Res}\nolimits}
\def\rank{\operatorname{rank}\nolimits}
\def\Cl{\operatorname{Cl}\nolimits}
\def\Pic{\operatorname{Pic}\nolimits}
\def\Picent{\operatorname{Picent}\nolimits}
\def\TrPic{\operatorname{TrPic}\nolimits}
\def\TrI{\operatorname{TrI}\nolimits}
\def\StPic{\operatorname{StPic}\nolimits}
\def\Stab{\operatorname{Stab}\nolimits}
\def\TrPicent{\operatorname{TrPicent}\nolimits}
\def\Shift{\operatorname{Shift}\nolimits}
\def\Grp{\operatorname{Grp}\nolimits}
\def\Spec{\operatorname{Spec}\nolimits}
\def\coker{\operatorname{coker}\nolimits}
\def\ie{{\em i.e.}}
\def\Rarr#1{\buildrel #1\over \longrightarrow}
\def\ov#1{\overline{#1}}
\def\iso{\buildrel \sim\over\to}
\def\Sh{\operatorname{Sh}\nolimits}
\def\Qlbar{{\bar{\bold{Q}}_l}}
\def\ti{{\tilde{i}}}
\def\ts{{\tilde{s}}}
\def\tB{{\widetilde{B}}}
\def\tH{{\widetilde{H}}}
\def\tI{{\tilde{I}}}
\def\tL{{\widetilde{L}}}
\def\tM{{\widetilde{M}}}
\def\tN{{\widetilde{N}}}
\def\tV{{\widetilde{V}}}
\def\tW{{\widetilde{W}}}
\def\tS{{\widetilde{S}}}
\def\tC{{\widetilde{C}}}
\def\tK{{\widetilde{K}}}
\def\tX{{\widetilde{X}}}
\def\tCP{{\widetilde{\CP}}}
\def\CExt{{{\mathcal E}xt}}
\def\CHom{{{\mathcal H}om}}
\def\SGt{{\tilde{\Sn}}}
\def\AGt{{\tilde{\GA}}}
\def\SG{{\Sn}}
\def\AG{{\GA}}

\def\DS{\displaystyle}

\def\a{\alpha}
\def\b{\beta}
\def\ga{\gamma}
\def\Ga{\Gamma}
\def\g{\gamma}
\def\G{\Gamma}
\def\d{\delta}
\def\D{\Delta}
\def\e{\varepsilon}
\def\ph{\varphi}
\def\ch{\chi}
\def\l{\lambda}
\def\L{\Lambda}
\def\m{\mu}
\def\n{\nu}
\def\o{\omega}
\def\O{\Omega}
\def\r{\rho}
\def\s{\sigma}

\def\th{\theta}
\def\Th{\Theta}
\def\t{\tau}
\def\x{\xi}
\def\X{\Xi}
\def\z{\zeta}

\def\matrice#1{\left(\begin{array}{ccccccccccccccccccc}#1\end{array}\right)}
\def\determinant#1{\left|\begin{array}{ccccccccccccccccccc}#1\end{array}\right|}
\def\tete#1{\par\leavevmode\makebox[0.7cm]{$(\mathrm{#1})$}}

\title{Quotients et extensions de groupes de r\'eflexion}
\author{David Bessis, C\'edric Bonnaf\'e et Rapha\"el Rouquier}

\begin{document}

\maketitle

\begin{quotation}
{\flushleft \bf Abstract.} {\small We give a geometric description
of a certain class of epimorphisms between complex reflection groups.
We classify these epimorphisms, which can be interpreted as
``morphisms'' between the diagrams symbolizing
standard presentations by generators and relations for complex
reflection groups and their braid groups.
}
\end{quotation}

\section{Introduction}

Dans l'\'etude des groupes de r\'eflexion, il est classique de consid\'erer
certaines classes de bons sous-groupes (les sous-groupes paraboliques) et
de bons automorphismes (les automorphismes ``de diagramme''). Nous
proc\'edons ici \`a la d\'efinition et \`a l'\'etude d'une classe de
bons quotients.

Notre motivation originale \'etait de comprendre pourquoi les
pr\'esentations par g\'en\'erateurs et relations classiques (voir les
tables de \cite{BrMaRou}) associ\'ees \`a diff\'erents groupes de
r\'eflexions pr\'esentent certaines similarit\'es.
Prenons l'exemple du diagramme de Coxeter de type $F_4$. Quand
on en supprime la double-barre, on obtient le diagramme de Coxeter
de type $A_2\times A_2$. En termes de pr\'esentations, supprimer la
double-barre consiste \`a ajouter une relation de commutation. Ainsi
cette op\'eration d\'ecrit un morphisme surjectif de
$W(F_4)$ vers $W(A_2\times A_2)$:
$$\xy (0,0) *++={\phantom{x}} *\frm{o} ; (10,0) *++={\phantom{x}} *\frm{o} **@{-} ;
 (20,0) *++={\phantom{x}} *\frm{o} **@{=}; (30,0) *++={\phantom{x}} *\frm{o} **@{-} ;
   (0,-3) *++={} ; (10,-3) *++={} ; (20,-3) *++={};
    (30,-3) *++={} \endxy \twoheadrightarrow
\xy (0,0) *++={\phantom{x}} *\frm{o} ; (10,0) *++={\phantom{x}} *\frm{o} **@{-} ;
 (20,0) *++={\phantom{x}} *\frm{o} ; (30,0) *++={\phantom{x}} *\frm{o} **@{-} ;
   (0,-3) *++={} ; (10,-3) *++={} ; (20,-3) *++={};
    (30,-3) *++={} \endxy$$
Quelle est la signification g\'eom\'etrique d'un tel morphisme~? Quelles
en sont les cons\'equences, en termes d'arrangements d'hyperplans,
d'invariants polynomiaux et de groupes de tresses ?

Des \'epimorphismes semblables
apparaissent fr\'equemment entre groupes de r\'eflexion complexes.
Nous expliquons comment construire des suites exactes
$$\xymatrix@1{ 1 \ar[r] & G \ar[r] & \tW \ar[r] & W \ar[r] & 1}$$
o\`u $\tW$ et $W$ sont deux groupes de r\'eflexion complexes.

\smallskip

Soit $\tW$ un groupe de r\'eflexion complexe, agissant sur
un espace $\tV$.
Dans la situation qui nous int\'eresse,
$G$ est un sous-groupe distingu\'e de $\tW$
et l'action de $W=\tW/G$ sur la vari\'et\'e quotient
$\tV/G$ s'\'etend
en une action lin\'eaire sur l'espace tangent $V$ \`a $\tV/G$ en
$0$.

Nous montrons que  $W$ est un groupe
de r\'eflexion sur $V$ si et seulement si $\tV/G$ est une
intersection compl\`ete et $W$ agit trivialement sur
$\Tor_*^{\BC[V]}(\BC[\tV/G],\BC)$.

Nous expliquons alors comment
relier les degr\'es de $\tW$, son arrangement
d'hyperplans, ses r\'eflexions, ses sous-groupes paraboliques, son groupe
de tresses
\`a ceux de $W$.
Par exemple, dans le cas (crucial) o\`u $G$ ne contient pas de r\'eflexions,
l'image d'une r\'eflexion de $\tW$ est une r\'eflexion de $W$ du m\^eme
ordre.

\smallskip
Nous abordons ensuite (partie 4)
le probl\`eme de la classification des paires $(\tW,G)$.
Un des r\'esultats obtenus est le suivant: 
lorsque $\tV$ est de dimension $2$ et $G=\{\pm 1\}$, la correspondance
$\tW\mapsto W$ induit une bijection entre les classes de conjugaison
de groupes de r\'eflexion complexes de dimension $2$ engendr\'es par
des r\'eflexions d'ordre $2$ et contenant $G$ et les groupes
de Coxeter finis (non triviaux) de dimension au plus $3$.

Nous donnons des tables de transformations \'el\'ementaires entre
diagrammes ``\`a la Coxeter'' qui permettent de d\'ecrire 
r\'ecursivement l'ensemble des paires $(\tW,G)$. 
Nous expliquons comment le morphisme $\tW \twoheadrightarrow W$ peut \^etre
interpr\'et\'e comme un ``morphisme de diagramme''.

\smallskip

Dans une derni\`ere partie, 
nous traitons le probl\`eme inverse~: au lieu de construire
$W$ \`a partir de $\tW$ et $G$, nous \'etudions
l'existence d'un groupe de r\'eflexion $\tW$
\'etant donn\'es $W$ et $G$.

\smallskip
La classification des groupes de r\'eflexion complexes
rec\`ele de nombreuses co\"\i ncidences num\'e\-rologiques et
correspondances internes. Certaines de
ces co\"\i ncidences sont bien comprises. Par exemple, certains groupes de
r\'eflexion sont des sous-groupes {\em paraboliques} d'autres groupes de
r\'eflexion. D'autres co\"\i ncidences sont expliqu\'ees par la th\'eorie
des { \'el\'ements r\'eguliers} de Springer (\cite{springer}). D'autres
encore peuvent \^etre intepr\'et\'ees gr\^ace \`a la {correspondance
de McKay}. Le pr\'esent travail permet d'expliquer un nouveau type de
correspondance entre groupes de r\'eflexion.

Outre le fait qu'elle permet de mieux comprendre la classification
de Shephard-Todd,
un des int\'er\^ets de la correspondance d\'ecrite ici est
de relier des groupes de Coxeter, dont la combinatoire est
bien connue, \`a certains groupes de r\'eflexions complexes (non r\'eels)
pour lesquels des pr\'esentations ne sont connues que de fa\c con empirique.

Nous d\'etaillons par exemple le cas du groupe exceptionnel $G_{31}$
qui admet comme quotient
de r\'eflexion le groupe sym\'etrique ${\mathfrak{S}}_6$, faisant
appara\^{\i}tre une hypersurface $\Sn_6$-invariante remarquable de $\BC^5$.
Or, si un diagramme pour $G_{31}$ est propos\'ee dans \cite{BrMaRou},
la question de savoir si ce diagramme donne aussi une pr\'esentation
du groupe de tresses associ\'e est toujours ouverte.
Ce diagramme est compatible, via le morphisme $G_{31} \twoheadrightarrow
{\mathfrak{S}}_6$, avec la pr\'esentation de Coxeter de ${\mathfrak{S}}_6$,
indice qui renforce la conviction qu'il est probablement le ``bon'' diagramme
pour $G_{31}$.

Reste bien entendu \`a pr\'eciser quelles propri\'et\'es 
exactes doivent v\'erifier de ``bons'' diagrammes des groupes de
r\'eflexion complexes. Le pr\'esent travail fournit
une nouvelle condition de compatibilit\'e \`a inclure dans le
cahier des charges.
 
\bigskip

\section{Pr\'eliminaires}~

\medskip

Soit $k$ un corps de caract\'eristique nulle. Si $n$ est un entier naturel non nul, 
on notera $\m_n(k)$ le groupe des racines $n^{\text{\`emes}}$ de l'unit\'e dans $k$. 

\bigskip

\subsection{}
Soit $W$ un groupe fini agissant sur une $k$-alg\`ebre $A$. Si $I$ est
un id\'eal de $A$, on notera $W_I$ (ou $W_{\Spec A/I}$)
le groupe d'inertie de $I$. Si $\chi$ est un caract\`ere
irr\'eductible de $W$, on notera $A^W_\chi$ la composante
$\chi$-isotypique de $A$.

\bigskip

\subsection{} Soit $V$ un espace vectoriel (tous les espaces vectoriels consid\'er\'es
seront de dimension finie) sur $k$. On notera $V^*$ le dual de $V$ et $k[V]=S(V^*)$, l'alg\`ebre 
sym\'etrique de $V^*$. Nous confondrons parfois le sch\'ema
$\Spec k[V]$ et l'espace vectoriel de ses points sur $k$, {\em i.e.},
$V$. Soit $I$ un id\'eal de $k[V]$. 
On dira que $I$ (ou $\Spec k[V]/I$) est une intersection compl\`ete si le nombre 
minimal de g\'en\'erateurs de $I$ est \'egal \`a $\dim V-\dim( \Spec k[V]/I)$.

Une graduation sur $V$ est une action de $k^\times$ 
avec poids strictement positifs. 
On dira qu'une graduation est standard si les poids
sont tous $1$. Si $V$ est muni d'une graduation, alors l'alg\`ebre $k[V]$ 
h\'erite d'une graduation, c'est-\`a-dire, d'une action de $k^\times$ avec poids
positifs ou nuls.

\bigskip

\subsection{} Soit $A$ une $k$-alg\`ebre commutative de type fini, munie d'une
graduation. Soit $A_n$ la composante de degr\'e $n$ de $A$. On pose
$A_+=\oplus_{i > 0} A_i$. On suppose que $A_0 = k$~: en particulier, 
$A_+$ est un id\'eal maximal de $A$. 
Soit $X=\Spec A$ et soit $0$ le point de $X$ correspondant \`a 
l'id\'eal maximal $A_+$ de $A$. Si $I$ est un id\'eal homog\`ene de $A$, alors 
le nombre minimal de g\'en\'erateurs homog\`enes de l'id\'eal $I$ est 
$\dim I/A_+ I$ par le lemme de Nakayama. 

L'espace tangent $V=(A_+/(A_+)^2)^*$ 
au sch\'ema $X$ en $0$ est un $k$-espace vectoriel gradu\'e. Le lemme 
suivant est classique~:

\begin{lemme}\label{generation}
Soient $p_1$,\dots, $p_n$ des \'el\'ements homog\`enes 
de $A_+$. Alors $A$ est engendr\'ee 
par $p_1$,\dots, $p_n$ si et seulement si les images de $p_1$,\dots, $p_n$ dans 
$V^*=A_+/(A_+)^2$ engendrent $V^*$. En particulier, $\dim V$ est le nombre minimal 
de g\'en\'erateurs de $A$.
\end{lemme}


D'apr\`es le lemme \ref{generation}, une 
scission gradu\'ee $V^* \to A_+$ du morphisme surjectif $A_+ \to V^*$ 
induit un morphisme surjectif de $k$-alg\`ebres gradu\'ees $k[V] \twoheadrightarrow A$, 
c'est-\`a-dire une immersion ferm\'ee de $\Spec A$ dans son espace tangent en $0$. 
Si $A$ est une alg\`ebre de polyn\^omes, alors l'application construite est 
un isomorphisme.

\bigskip

\subsection{}
Soit $\tV$ un $k$-espace vectoriel gradu\'e et soit $G$ un sous-groupe fini 
de $GL_{grad}(\tV)$, o\`u $GL_{grad}(\tV)$ est le sous-groupe de $GL(\tV)$ 
des \'el\'ements commutant \`a l'action de $k^\times$. On note $N(G)$ le normalisateur 
de $G$ dans $GL_{grad}(\tV)$~: $N(G)$ est un groupe r\'eductif
(non n\'ecessairement connexe) dont la 
composante connexe de l'\'el\'ement neutre est le centralisateur 
de $G$ dans $GL_{grad}(\tV)$.

L'alg\`ebre $k[\tV]^G$ est une $k$-alg\`ebre commutative 
gradu\'ee de type finie et $k[\tV]_0^G=k$. Le groupe $N(G)$ agit de mani\`ere 
gradu\'ee sur la $k$-alg\`ebre $k[\tV]^G$ et, puisqu'il est r\'eductif et
que $k$ est de caract\'eristique $0$, on peut 
choisir une scission gradu\'ee $N(G)$-\'equivariante $V^* \to k[\tV]_+^G$ 
du morphisme $k[\tV]_+^G \twoheadrightarrow V^*$, o\`u $V=(k[\tV]_+^G/(k[\tV]_+^G)^2)^*$.
On identifiera $\tV/G$ avec le sch\'ema $\Spec k[\tV]^G$. Alors $V$ est l'espace tangent 
\`a $\tV/G$ en $0$ et la scission pr\'ec\'edente induit une immersion
ferm\'ee $N(G)$-\'equivariante $\tV/G \to V$.

On note $I$ le noyau du morphisme surjectif d'alg\`ebres
$k[V] \twoheadrightarrow k[\tV]^G$, 
c'est-\`a-dire l'id\'eal de d\'efinition du sch\'ema $\tV/G$ dans $V$. 
On dit que $G$ est d'{\em intersection compl\`ete} si l'id\'eal $I$ l'est ou, 
de mani\`ere \'equivalente, si $\dim I/k[V]_+ I=\dim V - \dim \tV$ 
(en effet, $\krull k[\tV]^G=\dim \tV$).

L'ensemble des {\it degr\'es} de $G$ est l'ensemble des poids de $k^\times$
agissant sur $V$ (pris avec
multiplicit\'es). C'est aussi l'ensemble des degr\'es d'un syst\`eme
minimal d'invariants fondamentaux homog\`enes de $G$ dans son action
sur $k[\tV]$. 
L'ensemble des {\it degr\'es des relations} de $G$ est l'ensemble des
poids de $k^\times$ sur $I/k[V]_+ I$. On notera $N(G,rel)$ le 
sous-groupe de $N(G)$ des \'el\'ements qui agissent trivialement 
sur $I/k[V]_+ I$.

On appelle {\em r\'eflexion} de $V$ un automorphisme $g$ d'ordre fini tel que
$\ker (g-1)$ est un hyperplan de $V$.
On dira que $G$  est de {\em r\'eflexion}
s'il est engendr\'e par ses r\'eflexions.
Le th\'eor\`eme de Shephard-Todd-Chevalley \cite[Th\'eor\`eme 7.2.1]{Ben} 
donne plusieurs caract\'erisations des groupes de r\'eflexion en termes de 
leur alg\`ebre d'invariants~:

\begin{thm}[Shephard-Todd-Chevalley]\label{s-t-c} Les assertions suivantes sont \'equivalentes~:
\begin{itemize}
\item[(i)] $G$ est de r\'eflexion,

\item[(ii)] $k[V]^G$ est une alg\`ebre de polyn\^omes,

\item[(iii)] $k[V]$ est un $(k[V]^G)[G]$-module libre de rang $1$.
\end{itemize}
\end{thm}

L'\'equivalence entre les assertions (i) et (ii) du th\'eor\`eme \ref{s-t-c} 
montre qu'un groupe de r\'eflexion est d'intersection compl\`ete.

\smallskip

On appelle {\em double r\'eflexion} de $V$ un endomorphisme $g$ d'ordre fini
tel que $\ker (g-1)$ est de codimension $2$ dans $V$.
Kac et Watanabe \cite[Th\'eor\`eme A]{KaWa} ont
montr\'e que, si $G$ est d'intersection compl\`ete, alors 
$G$ est engendr\'e par ses r\'eflexions et ses
doubles r\'eflexions (la r\'eciproque n'est pas vraie).

\bigskip

\subsection{} Si $V$ est un $k$-espace vectoriel, 
$W$ est un groupe de r\'eflexion dans $GL_{grad}(\tV)$ et $L$ est un sous-espace de $V$, 
alors le
groupe $W_L$ est de r\'eflexion sur $V$ \cite[Th\'eor\`eme 1.5]{Ste}.
On note $\CA(W)$ l'ensemble des hyperplans des r\'eflexions de $W$. Si $W'$ est un groupe de
r\'eflexion sur $V'$, on dira que $f:W\to W'$ est un morphisme
de groupes de r\'eflexion si c'est un morphisme de groupes et si
l'image d'une r\'eflexion de $W$ est $1$ ou une r\'eflexion de $W'$.

\bigskip

\section{Quotients de groupes de r\'eflexion}~

\medskip

Soient $\tW$ un groupe de r\'eflexion sur un $k$-espace vectoriel $\tV$,
$G$ un sous-groupe distingu\'e de $\tW$ et $W=\tW/G$. 
Soit $V$ l'espace tangent \`a $\tV/G$ en $0$. On
fixe un plongement $N(G)$-\'equivariant $\tV/G\to V$.
On pose $B=k[\tV]$, $A=k[V]$ et on note $I$ l'id\'eal
de $A$ d\'efinissant la sous-vari\'et\'e ferm\'ee $\tV/G$.
On a $B^\tW=(A/I)^W=A^W/I^W$.

\medskip

On note $\CA'(\tW)$ le sous-ensemble des \'el\'ements $H$ de $\CA(\tW)$
tels que $\tW_H\not= G_H$. On d\'efinit
$\overline{\CA'}(\tW)$ comme le quotient de $\CA'(\tW)$
par la relation d'\'equivalence qui
identifie $H$ et $H'$ si $\tW_H$ et $\tW_{H'}$ ont la m\^eme
image dans $W$. Pour $H\in\CA(\tW)$, on fixe
une forme lin\'eaire $\alpha_H\in B$ d\'efinissant $H$.
Pour $C\subseteq\CA(\tW)$, on pose
$\alpha_C=\prod_{H\in C}\alpha_H^{e_H}$ o\`u $e_H=|G_H|$.

\bigskip

Dans cette partie, nous allons donner une condition n\'ecessaire et suffisante pour 
que le groupe $W$ soit de r\'eflexion sur $V$ (cf. th\'eor\`eme \ref{ci}).
Lorsque cette 
condition est r\'ealis\'ee, nous comparons les propri\'et\'es des groupes 
$W$ et $\tW$ (sous-groupes paraboliques, 
ordres des r\'eflexions, arrangements d'hyperplans, degr\'es des invariants, 
groupes de tresses associ\'es...).

\bigskip

\subsection{Bons sous-groupes distingu\'es des groupes de r\'eflexion} Dans ce paragraphe, 
nous allons montrer que $W$ est de r\'eflexion sur $V$ si et seulement si $G$
est d'intersection
compl\`ete et $W$ agit trivialement sur $I/A_+ I$. Cette deuxi\`eme condition a plusieurs 
interpr\'etations qui sont donn\'ees par le lemme suivant~:

\begin{lemme}
\label{red}
Les propri\'et\'es suivantes sont \'equivalentes~:
\begin{itemize}
\item[(i)]
le sch\'ema $\tV/\tW\times_{V/W}V$ est r\'eduit ;
\item[(ii)]
l'application canonique $\tV/G\to \tV/\tW\times_{V/W}V$ est
un isomorphisme ;
\item[(iii)]
on a $I=AI^W$ ;
\item[(iv)]
$W$ agit trivialement sur $\Tor_1^A(A/I,A/A_+)\simeq I/A_+ I$.
\end{itemize}
\end{lemme}

\begin{proof}
On a $\tV/\tW\times_{V/W}V=\Spec A^W/I^W\otimes_{A^W}A=\Spec
A/(AI^W)$. D'autre part, l'application canonique $\tV/G \to \tV/\tW \times_{V/W} 
V$ est bijective sur les points ferm\'es. Il r\'esulte alors du Nullstellensatz 
que l'application canonique
$A/(AI^W)\to A/I$ est un isomorphisme si et seulement si
$A/(AI^W)$ est r\'eduit. En outre, c'est un isomorphisme si et seulement
si $I$ est engendr\'e par des \'el\'ements $W$-invariants.

On a $\Tor_1^A(A/I,A/A_+) \simeq I/A_+I$. Ainsi, $W$ agit trivialement sur
ce module si et seulement si $I=AI^W+A_+ I$. D'apr\`es le lemme de Nakayama,
$I=A I^W+A_+ I$ implique $I=A I^W$.
\end{proof}

\bigskip

Le r\'esultat suivant est dans la m\^eme veine que les \S 3.3.2 et 3.3.3
de \cite{Go}.

\begin{thm}
\label{ci}
Le groupe $W$ est de r\'eflexion sur $V$ si et seulement si 
$G$ est d'intersection compl\`ete et $W$ agit trivialement sur $I/A_+ I$.
\end{thm}

\begin{proof}
Supposons que $W$ agit trivialement sur $I/A_+I$ et que $A/I$ est une intersection compl\`ete. 
Alors, $I$ peut \^etre
engendr\'e par $r$ \'el\'ements $W$-invariants, o\`u
$r=\dim V-\dim \tV$ (Lemme \ref{red}, \'equivalence entre (iii) et (iv)). 
Par cons\'equent, l'id\'eal $I^W$ de $A^W$
est engendr\'e par $r$ \'el\'ements. Puisque
$A^W/I^W=B^{\tW}$ est une alg\`ebre de polyn\^omes en $\dim \tV$ variables,
$A^W$ est engendr\'e par $\dim \tV + r $ \'el\'ements.
Comme $\krull A^W= \dim V = \dim \tV +r$, on en d\'eduit que
$A^W$ est une alg\`ebre de polyn\^omes, donc que $W$ est de
r\'eflexion sur $V$.

Supposons $W$ de r\'eflexion sur $V$. Alors,
$A$ est un $A^W$-module libre de rang $|W|$,
donc $A^W/I^W\otimes_{A^W}A$ est un
$A^W/I^W=B^\tW$-module libre de rang $|W|$.
Puisque $B\simeq B^\tW [\tW]$ comme $B^\tW [W]$-modules, on a 
$B^G\simeq B^{\tW}[W]$ comme $B^\tW[W]$-modules. Par cons\'equent,
$A/I=B^G$ est un $B^\tW$-module libre de rang $|W|$. Finalement, la surjection
canonique $A^W/I^W\otimes_{A^W}A\to A/I$ est un isomorphisme, donc $W$ 
agit trivialement sur $I/A_+I$ par le Lemme \ref{red}, \'equivalence entre (ii) et (iv).

Puisque $A^W$ et $A^W/I^W$ sont des alg\`ebres de polyn\^omes, $I^W$
est engendr\'e par $r$ \'el\'ements. Par cons\'equent, $I$ est
engendr\'e par $r=\krull A-\krull A/I$ \'el\'ements,
donc $I$ est une intersection compl\`ete.
\end{proof}

\bigskip

Le th\'eor\`eme \ref{ci} nous am\`ene \`a la d\'efinition suivante~:

\begin{defi}
On dira que $G$ est un bon sous-groupe distingu\'e de $\tW$
(ou simplement que $G$ est bon dans $\tW$) si $G$ est d'intersection
compl\`ete et si $W$ agit trivialement sur $I/A_+I$.
\end{defi}

Si $G$ est bon dans $\tW$, 
on a une action triviale de $W$ sur sur les groupes $\Tor$ sup\'erieurs~:

\begin{prop}
Si $G$ est un bon sous-groupe distingu\'e de $\tW$, alors on a un
isomorphisme de $A^W[W]$-modules~:
$\Tor_*^A(A/I,A/A_+)\simeq \Tor_*^{A^W}(A^W/I^W,A^W/(A^W)_+)$.
\end{prop}

\begin{proof}
On a 
$$A^W/I^W\otimes_{A^W}^\BL A\simeq
A^W/I^W\otimes_{A^W} A\simeq A/I$$
car $A$ est libre comme $A^W$-module et $I=AI^W$ car $G$ est bon ($\otimes^\BL$ 
d\'esigne le foncteur d\'eriv\'e \`a gauche de $\otimes$).
Puisque
$$(A^W/I^W\otimes_{A^W}^\BL A)\otimes_A^\BL A/A_+\simeq
A^W/I^W\otimes_{A^W}^\BL (A\otimes_A^\BL A/A_+),$$
on a
$$A/I\otimes_A^\BL A/A_+\simeq A^W/I^W\otimes_{A^W}^\BL A^W/(A^W)_+.$$
\end{proof}


\begin{prop}
\label{Grefl}
Si $G$ est de r\'eflexion, alors c'est un bon sous-groupe
distingu\'e de $\tW$.
\end{prop}

\begin{proof}
En effet, on a alors un isomorphisme $\tV/G\iso V$ et le th\'eor\`eme
\ref{ci} apporte la conclusion.
\end{proof}

\medskip

\begin{rem}
\label{cipasbon}
Si $G$ est d'intersection compl\`ete, il n'est pas n\'ecessairement
bon dans $\tW$. Pour $\tW=\mu_2(\BC)\times \mu_4(\BC)$ agissant
sur $\BC^2$ et $G$ le sous-groupe engendr\'e par $(-1,-1)$, $G$ est
d'intersection compl\`ete mais n'est pas bon dans $\tW$.
\end{rem}

\begin{rem}
\label{norm}
On suppose ici que $G$ n'est pas de r\'eflexion.
Si l'action de $G$ sur $\tV$ est absolument irr\'eductible,
alors $N(G,rel)$ est fini.

Soit $d$ le pgcd des degr\'es des relations de $G$. Alors, 
$N(G,rel)\cap k^\times={\mathbf \mu}_d(k)$. Par cons\'equent, 
si $\tW$ est irr\'eductible et $G$ est bon dans $\tW$, alors l'ordre du
centre de $\tW$ divise $d$ (en effet, un groupe de r\'eflexion est
irr\'eductible si et seulement si il est absolument irr\'eductible).
\end{rem}

\begin{rem}
Supposons $G$ d'intersection compl\`ete.
Alors, $G$ est bon dans $\tW$ si et seulement si
$\tW$ est contenu dans $N(G,rel)$. D'autre part, si $\tV/G$ est une 
hypersurface ({\em i.e.}, $\Tor_1^A(A/I,A/A_+)$ est de dimension $1$) et $k$ est 
alg\'ebriquement clos, alors $N(G)=k^\times\cdot N(G,rel)$.
\end{rem}


\subsection{Propri\'et\'es du quotient} 
On suppose ici que $G$ est un bon sous-groupe distingu\'e de $\tW$.
On sait que $Z=\tV/G$ est une intersection compl\`ete.
Notons $q:\tW\to W$ la surjection canonique et
$p:\tV\to \tV/G=Z\hookrightarrow V$. Nous allons \'etudier dans ce paragraphe 
les propri\'et\'es des applications $p$ et $q$.

\bigskip

\subsubsection{Sous-groupes paraboliques} La propri\'et\'e de la paire $(\tW,G)$ 
est h\'erit\'ee par les sous-groupes paraboliques~:

\begin{thm}
\label{parab1}
Soit $x\in \tV$. Alors, $G_x$ est bon dans $W_x$~: de plus, l'application canonique de 
l'espace tangent \`a $\tV/G_x$ en $x$ vers $V$ est injective et elle induit un isomorphisme
de groupes de r\'eflexion $\tW_x/G_x\iso W_{p(x)}$.
\end{thm}

\begin{proof}
Soit $\Gm$ l'id\'eal
maximal de $B$ d\'efinissant $x$. Notons $H=G_x$ le groupe de
d\'ecomposition de $\Gm$ dans $G$,
$\tK=\tW_x$ le groupe de d\'ecomposition
de $\Gm$ dans $\tW$ et $K=\tK/H$.

Soit $\Gm_H=\Gm\cap B^H$ et $\Gm_G=\Gm\cap B^G$. Notons que
$K$ est le groupe de d\'ecomposition de $\Gm_G$.

Le morphisme
$\Spec B^H_{\Gm_H}\to \Spec B^G_{\Gm_G}$ est \'etale, donc il
induit un isomorphisme de $kK$-modules
$(\Gm_H/\Gm_H^2)^*\iso(\Gm_G/\Gm_G^2)^*$. Il s'agit donc de d\'emontrer
que $K$ est de r\'eflexion sur $E=(\Gm_G/\Gm_G^2)^*$.

Soit $\Gn$ l'id\'eal maximal de
$A$ au-dessus de $\Gm_G$.
L'application naturelle de l'espace tangent
$E$ de $A_\Gn/I_\Gn$ vers l'espace tangent de $A_\Gn$
est injective. Cette application est $K$-\'equivariante.
L'espace tangent de $A_\Gn$ est canoniquement isomorphe \`a $V$ (et cet
isomorphisme est compatible \`a l'action de $K$).
Puisque $K$ agit fid\`element sur
$E$ et est de r\'eflexion sur $V$ (c'est le groupe de d\'ecomposition
de $\Gn$), il est de r\'eflexion sur $E$, d'o\`u le r\'esultat.
\end{proof}

\begin{cor}
\label{parab2}
Soit $L$ un sous-espace de $\tV$.
Alors, $G_L$ est un bon sous-groupe distingu\'e de $W_L$.
L'application canonique $\tW_L/G_L\to W_{p(L)}$ est un 
isomorphisme de groupes de r\'eflexion.
\end{cor}

\begin{proof}
Choisissons un point $x\in L$ tel que $\tW_x=\tW_L$. 
Les espaces tangents \`a $V/G_x$ en $x$ et $0$ sont isomorphes de mani\`ere
$\tW_x$-\'equivariante, donc $\tW_x/G_x$ est un groupe
de r\'eflexion sur l'espace tangent \`a $\tV/G_x$ en $0$, d'apr\`es le
th\'eor\`eme \ref{parab1} et l'application canonique
$\tW_x/G_x\to W_{p(x)}$ est un isomorphisme qui envoie une r\'eflexion
sur une r\'eflexion.
\end{proof}

\bigskip

\begin{rem}
Si on prend $\tW'$ un sous-groupe de r\'eflexion quelconque de $\tW$, alors
$G\cap\tW'$ n'est pas n\'ecessairement bon dans $\tW'$ (prendre
$\tW=\mu_4(\BC)\times\mu_4(\BC)$, $\tW'=\mu_4(\BC)\times\mu_2(\BC)$ et
$G=\tW\cap SL_2(\BC)$). 
\end{rem}


\bigskip

\subsubsection{Comparaison entre $W$ et $\tW$}

\begin{thm}
\label{prop}
On a les propri\'et\'es suivantes~:
\begin{itemize}
\item[(i)]
 $G$ est engendr\'e par des r\'eflexions et des doubles r\'eflexions.
\item[(ii)]
 Le morphisme $\tW\to W$ est un morphisme de groupes de r\'eflexion.
 L'image d'une r\'eflexion $s$ d'ordre $r$ autour d'un hyperplan $\tH$ est 
 une r\'eflexion d'ordre $r/(r\wedge e_\tH)$ si $r\not|\ e_\tH$ ; son
d\'eterminant est alors $\det(s)^{e_\tH}$.
\item[(iii)]
 Soit $\tH\in\CA'(\tW)$. Alors,
 il existe un unique hyperplan de r\'eflexion $\phi(\tH)$ de $W$
 contenant $p(\tH)$. L'application $\phi$ est un isomorphisme
 $\overline{\CA'}(\tW)\iso \CA(W)$. Pour $H\in\CA(W)$, l'ensemble $\phi^{-1}(H)$
 est la r\'eunion des composantes irr\'eductibles de
 $p^{-1}(H\cap Z)$.
\item[(iv)]
 Pour $H\in\CA(W)$, il existe une forme lin\'eaire $\alpha_H\in A$
 d\'efinissant $H$ dont l'image dans $A/I$ est $\alpha_{\phi^{-1}(H)}$.
 L'alg\`ebre $B^G$ est engendr\'ee par
 $\{\alpha_{\phi^{-1}(H)}\}_{H\in\CA(W)}$ et $B^\tW$.
\item[(v)]
 Soit $H$ un hyperplan de r\'eflexion de $W$ et $H_1$, $H_2$ deux
 composantes irr\'eductibles distinctes de $p^{-1}(H\cap Z)$. Alors,
 $H_1\cap H_2$ est contenu dans le lieu de ramification de $G$.
\end{itemize}
\end{thm}

\begin{proof}
Puisque $\tV/G$ est une intersection compl\`ete, (i) r\'esulte de
\cite[Th\'eor\`eme A]{KaWa}.
\smallskip

L'assertion (ii) est une cons\'equence imm\'ediate du corollaire
\ref{parab2}.

\smallskip

Soit $\tH\in\CA'(\tW)$. D'apr\`es (ii) et le corollaire \ref{parab2},
il existe un unique hyperplan de r\'eflexion $\phi(\tH)$ contenant $p(\tH)$.
Puisque l'image de $\tH$ dans $\overline{\CA}(\tW)$ est 
d\'etermin\'ee par $W_{\phi(\tH)}$, l'application $\phi$ est
injective.

Pour $H\in\CA(W)$, $H\cap Z$ est purement de codimension $1$ dans $Z$.
La sous-vari\'et\'e $p^{-1}(H\cap Z)$ de $\tV$ est une r\'eunion (non vide)
d'hyperplans de r\'eflexion, car c'est une sous-vari\'et\'e
ferm\'ee purement de codimension $1$ de $\tV$ contenue dans le lieu de
ramification de $\tW$. Cela prouve la surjectivit\'e de $\phi$ et donc
(iii).

\smallskip

Soit $H\in\CA(W)$ et $\alpha_H\in A$ une forme lin\'eaire d\'efinissant $H$.
Alors, l'image $\alpha$ de $\alpha_H$ dans $A/I$ est, \`a une constante pr\`es,
un produit $\prod_{\tH\in\phi^{-1}(H)}\alpha_\tH^{a_\tH}$ avec $a_\tH\ge 1$.
Soit $\tH\in\phi^{-1}(H)$, $x\in\tH$ tel que $\tW_\tH=\tW_x$ et
$\Gm$ l'id\'eal maximal correspondant de $B$. Soit $\Gm_G=\Gm\cap B^G$ et
$\Gm_K=\Gm\cap B^K$ o\`u $K=G_x$.
Soit $\Gn$ l'id\'eal maximal de $A$ au dessus de $\Gm_G$. 
Le noyau de l'application canonique $\Gn/\Gn^2\to \Gm_G/\Gm_G^2$ est
contenu dans les invariants par $\tW_x$ (car $\tW_x/K$ est un groupe
de r\'eflexion sur $\Gn/\Gn^2$ et agit fid\`element sur $\Gm_G/\Gm_G^2$).
Par cons\'equent, l'image de
$\alpha_H+\Gn^2$ par cette application est non-nulle. Ainsi,
l'id\'eal $\alpha B^G_{\Gm_G}$ de l'anneau local r\'egulier $B^G_{\Gm_G}$
est premier.
Puisque $B^K_{\Gm_K}$ est \'etale sur $B^G_{\Gm_G}$, l'id\'eal
premier $\alpha_\tH^{e_\tH}B^K_{\Gm_K}$ de $B^K_{\Gm_K}$ est engendr\'e par
$\alpha$.
Par cons\'equent, $\alpha B_\Gm=\alpha_\tH^{e_\tH} B_\Gm$, donc $a_\tH=e_\tH$.

Puisque $\{\alpha_H\}_{H\in\CA(W)}$ et $A^W$ engendrent $A$,
on a \'etabli (iv). 

\smallskip

Soient $s_1$ et $s_2$ deux r\'eflexions de $\tW$ d'hyperplans
de r\'eflexion $H_1$ et $H_2$ telles que $q(s_1)=q(s_2)$ ---
cela existe par (ii). Soit $g=s_1s_2^{-1}$. On a $g\in G$ et $g$
n'est pas trivial, puisqu'il agit non trivialement sur
$H_1\setminus(H_1\cap H_2)$.
Comme $g$ agit trivialement sur $H_1\cap H_2$, on en d\'eduit (v).
\end{proof}

\medskip

\begin{rem}
La sous-vari\'et\'e $p(\tH)$ de $V$ peut \^etre contenue
dans un sous-espace vectoriel strict de $\phi(\tH)$ (prendre
$\tW=\mu_2(\BC)\times\mu_2(\BC)$ et $G$ engendr\'e par $(-1,-1)$).
En outre, la vari\'et\'e $H\cap Z$ n'est pas n\'ecessairement
irr\'eductible (m\^eme exemple).
\end{rem}

\smallskip

\begin{cor}
\label{propstrict}
Supposons que $G$ ne contient pas de r\'eflexion (c'est-\`a-dire,
$e_\tH=1$ pour $H\in\CA(\tW)$). Alors,
\begin{itemize}
\item[(i)]
$G$ est engendr\'e par des doubles r\'eflexions.
\item[(ii)]
 Une r\'eflexion de $\tW$
 a pour image dans $W$ une r\'eflexion du m\^eme ordre.
\item[(iii)]
$G$ est contenu dans $SL(\tV)$.
\item[(iv)]
 L'image de l'id\'eal Jacobien de $A^W$ sur $A$ est l'id\'eal Jacobien de
 $B^\tW$ sur $B$.
\end{itemize}
\end{cor}

\begin{proof}
Les assertions (i) et (ii) sont des cons\'equences imm\'ediates
du th\'eor\`eme \ref{prop}.

Puisque $\alpha_C\in B^G$ pour $G\in\CA(\tW)$ (th\'eor\`eme \ref{prop} (iv)),
on a $\alpha_{\CA(\tW)}\in B^G$. Or,
$\alpha_{\CA(\tW)}\in B^{\tW}_{\det^{-1}}$ (cf. par exemple le lemme
\ref{calcul}), donc $\det_{|G}=1$, ce qui montre (iii).

L'assertion (iv) est une cons\'equence du th\'eor\`eme \ref{prop} (iv), car
l'id\'eal Jacobien de $A^W$ sur $A$ est engendr\'e par 
$\prod_{H\in\CA(W)} \alpha_H^{|W_H|-1}$ et
l'id\'eal Jacobien de $B^\tW$ sur $B$ par
$\prod_{\tH\in\CA(\tW)} \alpha_\tH^{|\tW_\tH|-1}$.
\end{proof}

\medskip

\begin{rem}
Soit $H$ le sous-groupe de $G$ engendr\'e par les r\'eflexions de $G$.
Quitte \`a remplacer $\tW$, $G$ et $\tV$ par $\tW/H$,
$G/H$ et $\tV/H$, on se ram\`ene gr\^ace \`a la proposition
\ref{Grefl} au cas o\`u
$G$ ne contient pas de r\'eflexion, c'est-\`a-dire, au cas o\`u l'application 
$\tV\to \tV/G$ est non ramifi\'ee en codimension $1$.
\end{rem}

La relation entre les degr\'es de $\tW$, $W$ et $G$ est donn\'ee par la
proposition imm\'ediate suivante (cf. remarque \ref{norm})~:

\begin{prop}\label{degres}
Munissons $\tV$ de sa graduation standard et $V$ de la graduation induite.
Soit $V_i$ le sous-espace de $V$ de degr\'e $i$, $W_i$ le fixateur
de $\oplus_{j\not=i}V_j$ et $E_i$ l'ensemble des degr\'es
de $W_i$ pour la graduation standard sur $V_i$.

Alors, $W=\times_i W_i$ et $W_i$ est un groupe
de r\'eflexion sur $V_i$. En outre,
$\bigcup_i iE_i$ est la r\'eunion de l'ensemble des degr\'es de $\tW$
avec l'ensemble des degr\'es des relations de $G$.
\end{prop}

\medskip

Dans le cas d'un groupe de Coxeter $\tW$, le groupe $W$ est naturellement
un groupe de Coxeter~:

\begin{prop}
\label{coxeter}
Supposons que $(\tW,\tS)$ est un groupe de Coxeter dans sa repr\'esentation
naturelle. Soit
$S$ l'ensemble des \'el\'ements non triviaux de l'image de $\tS$ dans $W$.
Alors, $(W,S)$ est un groupe de Coxeter.
\end{prop}

\begin{proof}
Nous avons ici $k=\BR$. Dans cette preuve, nous consid\'ererons la topologie
classique. 
Soit $X$ (resp. $\tX$) le compl\'ementaire des hyperplans de r\'eflexion de
$W$ (resp. $\tW$) dans $V$ (resp. $\tV$).
Alors, $\tS$ est l'ensemble des r\'eflexions par rapport
aux murs d'une chambre $\tC$ (=composante connexe) de $\tX$. Soit
$D$ l'image de $\tC$ dans $V$. Alors, $D$ est connexe et on note 
$C$ la composante connexe de $X$ la contenant.

Soit $\tH$ un hyperplan
de r\'eflexion d'un \'el\'ement $\ts\in\tS$ tel que $\ts\notin G$.
Soit $\tH^+$ l'intersection de $\tH$ avec l'adh\'erence de $\tC$.
Alors, $p(\tH^+)$ est contenu dans l'adh\'erence de $C$. Puisque
$p(\tH^+)$ est contenu dans un unique hyperplan de r\'eflexion
$\phi(\tH)$ de $W$ (th\'eor\`eme \ref{prop} (iii)), cet hyperplan est un
mur de la chambre $C$. 
Par cons\'equent, l'image de $\ts$ est une r\'eflexion autour d'un mur
de $C$.
\end{proof}

\bigskip

\subsubsection{Groupes de tresses} Dans ce paragraphe, nous supposons que
$k=\BC$. La topologie consid\'er\'ee sera la topologique classique. Posons~: 
$$\widetilde{\CM}=\tV-\bigcup_{\tH\in\CA(\tW)} \tH$$ $$\CM=V-\bigcup_{H\in\CA(W)} H.\leqno{\mathrm{et}}$$ 
Soit $x_0\in \widetilde{\CM}$, $\tB=\pi_1(\widetilde{\CM}/\tW,x_0)$ et
$B=\pi_1(\CM/W,p(x_0))$. D'apr\`es le th\'eor\`eme \ref{prop}, (iii), 
l'application $p : \tV \to V$ v\'erifie $p(\widetilde{\CM})\subset \CM$ donc elle 
induit un morphisme de groupes $p_*:\tB\to B$. 

Soit $\tH\in\CA(\tW)$, $\gamma$ un chemin de $\tV$ tel que
$\gamma(0)=x_0$, $\gamma(t)\in\tilde{\CM}$ pour $t<1$ et
$\gamma(1)\in \tH-\bigcup_{\tH'\in\CA(\tW),\tH'\not=\tH}\tH\cap \tH'$.
Soit $s_\tH$ le g\'en\'erateur de d\'eterminant
$\exp(2i\pi/|\tW_\tH|)$ de $\tW_\tH$.
On note $[\gamma]\in\tB$ la classe d'homotopie de  chemins de $x_0$ \`a
$s_\tH(x_0)$ dans $\tilde{\CM}$ associ\'ee \`a $\gamma$
\cite[Appendix 1 et \S 2.13]{BrMaRou}~: c'est un $s_\tH$-g\'en\'erateur
de la monodromie.

Alors $p([\gamma])=[p(\gamma)]$. Si $\tH\in \CA'(\tW)$,
c'est un $q(s_\tH)$-g\'en\'erateur de la monodromie.
Sinon, c'est l'\'el\'ement neutre de $B$.

\smallskip
Supposons maintenant que $\tW$ est
le complexifi\'e d'un groupe de r\'eflexion r\'eel
agissant sur l'espace vectoriel r\'eel $\tV'$, avec
$\tV=\tV'\otimes_\BR \BC$.
On reprend les notations de la proposition \ref{coxeter} et de sa preuve.
On fixe une chambre $\tC$ de $\tW$ dans $\tV'$ et on choisit
un point $x_0$ de $\tC$.
Pour $\ts\in\tS$, on choisit un chemin $\gamma$ dans $\tV'$ avec
$\gamma(0)=x_0$, $\gamma(t)\in C$ pour $t<1$ et $\tW_{\gamma(1)}=<\ts>$.
Soit $\sigma_{\ts}=[\gamma]$ l'\'el\'ement de $\tB$ associ\'e.
Brieskorn \cite{Br} a montr\'e (voir aussi Deligne \cite[Th\'eor\`eme 4.4]{De})
que $\tB$ est engendr\'e par les $\sigma_{\ts}$, avec $\ts\in \tS$ et que les
relations entre les $\sigma_{\ts}$ sont les relations de tresses.

De m\^eme, on a des g\'en\'erateurs $\sigma_s$ de $B$ pour
$s\in S$~: d'apr\`es ce qui pr\'ec\`ede, ce sont les images des 
$\sigma_{\ts}$, pour $\ts\in\tS$, $\ts\notin G$.

On a montr\'e le r\'esultat suivant~:

\begin{prop}
Si $\tW$ est le complexifi\'e d'un groupe de r\'eflexion r\'eel, alors
l'application canonique $p_*:\tB\to B$ est une surjection. Elle
envoie $\sigma_{\ts}$ sur $\sigma_{q(\ts)}$ si $\ts\notin G$
et sur $1$ sinon.
\end{prop}

Cette proposition reste vraie dans le cas g\'en\'eral o\`u
$\tW$ n'est pas n\'ecessairement le complexifi\'e d'un groupe
de r\'eflexion r\'eel
(voir \cite{Bes}).

\subsubsection{Extensions successives et r\'eductibilit\'e}
La proposition suivante est imm\'ediate.

\begin{prop}\label{prop sous-groupe}
Soit $\tW'$ un sous-groupe de r\'eflexion de $\tW$ contenant $G$.
Alors, $G$ est un bon sous-groupe distingu\'e de $\tW'$.
\end{prop}

\begin{rem}
Il se peut que $G$ soit un bon sous-groupe distingu\'e de $\tW'$ sans
\^etre un bon sous-groupe distingu\'e de $\tW$ (prendre
$\tW'=\mu_2(\BC)\times\mu_2(\BC)$ dans la remarque \ref{cipasbon}).
\end{rem}

\begin{prop}
\label{quotient}
Soit $H$ un sous-groupe distingu\'e de $\tW$ contenant $G$.
Alors, $H$ est bon dans $\tW$ si et seulement si 
$H/G$ est bon dans $W$.

Si $H$ est bon dans $\tW$, l'isomorphisme canonique
$\tW/H\iso W/(H/G)$ est un isomorphisme de groupes de r\'eflexion.
\end{prop}

\begin{proof}
Soit $\Gm$ l'id\'eal maximal de $A^H$ d\'efinissant le point $0$ de $V/H$.
L'injection de $\tV/G$ dans $V$ induit une injection de 
du plan tangent $(\Gm/(I^H+\Gm^2))^*$ \`a $\tV/H$ en $0$
dans le plan tangent $(\Gm/\Gm^2)^*$ \`a $V/H$ en $0$.
La proposition d\'ecoulera de la propri\'et\'e de $\tW$ d'agir
trivialement sur le conoyau de cette injection.

On a $I=I^{\tW}A$ (th\'eor\`eme \ref{ci}),
donc $I^H=I^{\tW}A^H=I^\tW+I^\tW \Gm\subseteq
I^\tW+\Gm^2$ et finalement $I^H+\Gm^2=I^\tW+\Gm^2$. Par
cons\'equent, $\tW$ agit trivialement sur
$(I^H+\Gm^2)/\Gm^2$ qui est le noyau de la
surjection canonique $\Gm/\Gm^2\to \Gm/(I^H+\Gm^2)$.
\end{proof}

\bigskip

Les propositions \ref{quotient} et \ref{Grefl} ram\`enent la 
classification des paires $(\tW,G)$ (avec $G$ bon) au cas
o\`u $G$ est de r\'eflexion ou bien ne contient pas de r\'eflexion.

\smallskip

Si le groupe $\tW$ est r\'eductible, une d\'ecomposition de $\tW$ est presque
toujours compatible \`a une d\'ecomposition de $G$~:

\begin{prop}
\label{reduc}
Soient $\tV_1$ et $\tV_2$ des sous-espaces $\tW$-stables de $\tV$ avec
$\tV=\tV_1\oplus \tV_2$. Soient $\tW_1=\tW_{\tV_2}$ et $\tW_2=\tW_{\tV_1}$.
Alors $G_1=G\cap\tW_1$ (resp. $G_2=G\cap\tW_2$) est un bon sous-groupe
distingu\'e de $\tW_1$ (resp. $\tW_2$).  En outre,
$G/(G_1\times G_2)$ est un bon sous-groupe distingu\'e de
$\tW/(G_1\times G_2)$.

Si $G_1=G_2=1$, $G\not=1$ et $\tV_1$ et $\tV_2$ sont irr\'eductibles,
alors, $\tV$ est de dimension $2$, il existe
un entier $d$ tel que $\tW_1$ et $\tW_2$ sont
des groupes de r\'eflexion cycliques d'ordre $d$ et $G=\tW\cap SL(\tV)$.
\end{prop}

\begin{proof}
D'apr\`es le corollaire \ref{parab2}, $G_1$ et $G_2$ sont
bons dans $\tW_1$ et $\tW_2$. Par cons\'equent, $G_1\times G_2$ est bon
dans $\tW$. Il r\'esulte maintenant de la proposition \ref{quotient}
que $G/(G_1\times G_2)$ est bon dans $\tW/(G_1\times G_2)$.

Supposons maintenant $G_1=G_2=1$, $G\not=1$ et $\tV_1$, $\tV_2$
irr\'eductibles. Examinons tout d'abord le cas o\`u $\tV_1$ et $\tV_2$ 
sont de dimension $1$. Alors $\tW$ est ab\'elien et $G$ ne contient pas de 
r\'eflexion. Par cons\'equent, $G \subseteq SL(\tV)$ (cf. 
corollaire \ref{propstrict}) donc est cyclique. Soit 
$d$ son ordre. Alors $\tW_1$ et $\tW_2$ sont cycliques 
d'ordre $d$ comme le montre le calcul de $N(G,rel)$ 
qui sera effectu\'e aux paragraphes \ref{sub C2} et \ref{subsub Cd}. 
En particulier, $G=\tW \cap SL(\tV)$ et est engendr\'e 
par un \'el\'ement de la forme $(s_1,s_2) \in \tW_1 \times 
\tW_2 = \tW$ o\`u $s_1$ et $s_2$ sont d'ordre $d$.

Revenons maintenant au cas o\`u $\tV$ est de dimension quelconque. 
Soit $H_1$ un hyperplan de r\'eflexion de $\tW_1$ agissant
sur $\tV_1$, $H_2$ un hyperplan de r\'eflexion de $\tW_2$ agissant
sur $\tV_2$ et $L=H_1\oplus H_2$. Alors, $G_L$ est
un bon sous-groupe distingu\'e de $\tW_L=\tW_{H_1}\times \tW_{H_2}$ dans son
action sur $(\tV_1/H_1)\oplus (\tV_2/H_2)$ qui est de dimension $2$ (cf. corollaire 
\ref{parab2}).
D'apr\`es la remarque pr\'ec\'edente, on en d\'eduit que
$G_L$ est engendr\'e par un \'el\'ement $g=(s_1,s_2)$ o\`u
$s_i$ est une r\'eflexion g\'en\'eratrice de $\tW_{\tH_i}$.

Soit $w\in \tW_1$. Alors, $gwg^{-1}w^{-1}\in G\cap \tW_1=1$. Par
cons\'equent, $s_1\in Z(\tW_1)$. La repr\'esentation $\tV_1$ \'etant
irr\'eductible, elle est par cons\'equent de dimension $1$ et
$s_1$ engendre $\tW_1$. De m\^eme, $\tV_2$ est de dimension $1$ et
$\tW_2$ est engendr\'e par $s_2$.
\end{proof}

Les deux propositions pr\'ec\'edentes permettent de ramener la
classification des paires $(\tW,G)$ (avec $G$ bon) au cas o\`u $\tW$ est
irr\'eductible.

%
%
%
%

\subsection{Construction \`a partir du rang $2$}
\label{fromrang2}
Le r\'esultat suivant pr\'ecise le th\'eor\`eme \ref{prop} et
explique comment construire $G$ bon dans $\tW$.

\begin{thm}
\label{constr}
Les assertions suivantes sont \'equivalentes~:
\begin{itemize}
\item[(i)]
le groupe $G$ est bon dans $\tW$ ;
\item[(ii)]
il existe une
famille $\CG$ de sous-groupes de $G$ engendrant $G$ et pour chaque
$G'\in \CG$, un sous-espace $L$ de $\tV$ tel que $G'$ est un bon 
sous-groupe distingu\'e de $\tW_L$ ;
\item[(iii)]
il existe une
famille $\CG$ de sous-groupes de $G$ engendrant $G$ et pour chaque
$G'\in \CG$, un sous-espace $L$ de
codimension $2$ de $\tV$ tel que $G'$ est un bon 
sous-groupe distingu\'e de $\tW_L$.
\end{itemize}
\end{thm}

Commen\c{c}ons par d\'ecrire l'action d'une r\'eflexion ou d'une double r\'eflexion
sur un produit de formes lin\'eaires.

\begin{lemme}
\label{calcul}
Soit $\CF$ une famille finie d'hyperplans de $\tV$ stable et permut\'ee
transitivement par un \'el\'ement $g\in GL(\tV)$. Soit
$\alpha=\prod_{H\in\CF}\alpha_H$
o\`u $\alpha_H$ est une forme lin\'eaire d\'efinissant $H$.
Soit $H\in\CF$ et $a\in k^\times$ d\'efini
par $g^{|\CF|}(\alpha_H)=a\alpha_H$.
Alors, $g(\alpha)=a\alpha$.

\smallskip
Supposons $\CF$ r\'eduite \`a un hyperplan $H$ et $a\not=1$,
{\em i.e.}, $g(\alpha_H)\in (k-\{1\})\alpha_H$.

Si $g$ est une r\'eflexion, alors $H=\ker (g-1)$.

Si $g$ est une double r\'eflexion,
alors il existe une (unique) r\'eflexion $s$ d'hyperplan $H$
telle que $gs$ est une r\'eflexion. En particulier, $\ker (g-1)\subset H$.
\end{lemme}

\begin{proof}
Soit $a_H\in k^\times$ d\'efini par $g(\alpha_H)=a_H\alpha_{H'}$
o\`u $H'=g(H)$. Alors, $g(\alpha)=(\prod_{H'\in\CF}a_{H'})\alpha$ et
$g^{|\CF|}(\alpha_H)=(\prod_{H'\in\CF} a_{H'})\alpha_H$, d'o\`u la
premi\`ere partie du lemme.

L'hyperplan $H$ est stable par $g$ et
$g$ n'agit pas trivialement sur $\tV/H$.

Si $g$ est une r\'eflexion on en d\'eduit que $H=\ker (g-1)$.

Supposons maintenant que $g$ est une double r\'eflexion.
Il existe une unique r\'eflexion $s$ d'hyperplan $H$ telle que 
$gs(\alpha_H)=\alpha_H$. Alors, $gs$ est une r\'eflexion.
\end{proof}

Le lemme suivant fournit un crit\`ere pour que $G$ soit bon en termes 
d'invariants.

\begin{lemme}
\label{invar}
Supposons que $G$ ne contient pas de r\'eflexion. Alors, $G$ est bon dans
$\tW$ si et seulement si $G\subseteq SL(\tV)$ et
pour tout $C\in\overline{\CA'}(\tW)$, on a
$\alpha_C\in B^G$.
\end{lemme}

\begin{proof}
Supposons $G$ bon dans $\tW$. Les propri\'et\'es requises sont donn\'ees par
le th\'eor\`eme \ref{prop} (iv) et le corollaire \ref{propstrict} (iii).

\smallskip
Montrons maintenant la r\'eciproque. Soit $H\in\CA(\tW)$,
$C\in\overline{\CA'}(\tW)$ contenant $H$ et $K=G\tW_H$.
L'ensemble $C$ est stable sous l'action de $K$. 
D'apr\`es le lemme \ref{calcul}, $\alpha_C\in B^{\tW_H}_{\det^{-1}}$.
Donc, $\alpha_C\in B^K_{\det^{-1}}$.
Il r\'esulte de \cite[Lemme 2.2]{Sta} que la famille des hyperplans de
r\'eflexion de $K$ est incluse dans $C$. Puisque les \'el\'ements
de $C$ sont des hyperplans de r\'eflexion de $K$, on en d\'eduit
que $\CA(K)=C$.

Pour $s'$ une r\'eflexion autour d'un hyperplan $H'\in C$, il existe
$s\in \tW_H$ tel que $ss'\in G$, donc $\det(ss')=1$ et $s$, $s'$
ont le m\^eme ordre. Par cons\'equent, $|\tW_{H'}|\le |\tW_H|$. On
montre de m\^eme l'in\'egalit\'e inverse. Par cons\'equent,
l'entier $|\tW_{H'}|$ est ind\'ependent du choix de $H'\in C$. On le
note $e$.

Puisque $\alpha^i_C\in B^K_{\det^{-i}}$, il r\'esulte de
\cite[Th\'eor\`eme 2.3]{Sta} que 
$B^K_{\det^{-i}}=\alpha^i_C B^K$ pour $0\le i<e$.

On a montr\'e que
$B^G=\bigoplus_{0\le i<e} \alpha^i_C B^K$, donc
$B^G$ est engendr\'ee par $B^K$ et $\alpha_C$.
Par cons\'equent, l'espace des invariants de $\tW_H$ dans $V$
est un hyperplan et donc $\tW_H$ est un groupe de r\'eflexion sur $V$.
Donc, $W$ est de r\'eflexion sur $V$.
\end{proof}

\begin{proof}[D\'emonstration du th\'eor\`eme \ref{constr}]
Notons que (iii)$\Rightarrow$(ii) est clair.

Si $G$ est bon, alors il est engendr\'e par des r\'eflexions et des
doubles r\'eflexions d'apr\`es le th\'eor\`eme \ref{prop}, donc
il existe une famille finie $\CF$ de sous-espaces de codimension $2$ de $\tV$
telle que $G$ est engendr\'e par les $G_L$ pour $L\in\CF$. D'apr\`es
le corollaire \ref{parab2}, $G_L$ est un bon sous-groupe distingu\'e de
$\tW_L$. On a montr\'e (i)$\Rightarrow$(iii).

\smallskip
Montrons maintenant (ii)$\Rightarrow$(i).
Gr\^ace aux propositions \ref{quotient} et
\ref{Grefl} et au corollaire \ref{parab2}, on peut remplacer $G$ et $\tW$ par
$G/G_r$ et $\tW/G_r$, o\`u $G_r$ est le sous-groupe engendr\'e par
les r\'eflexions de $G$. On peut ainsi supposer que $G$ ne contient
pas de r\'eflexion.

Soit $C\in\overline{\CA'}(\tW)$.
Soit $G'\in\CG$, $L$ un sous-espace de $\tV$ tel que $G'$ est bon dans
$\tW_L$ et $g\in G'$ une double r\'eflexion.

On a $\alpha_C=\alpha_{C'}\alpha_{C''}$
o\`u $C'=C\cap\CA(\tW_L)$ et $C''=C-C'$.
L'ensemble $C'$ admet une partition en parties
$D$ form\'ee des hyperplans $H'$ avec $G'\tW_{H'}$ fix\'e. Pour une
telle partie $D$, il r\'esulte du lemme \ref{invar} que $\alpha_D\in B^g$.
Puisque $\alpha_{C'}$ est un produit de tels $\alpha_D$, on a
$\alpha_{C'}\in B^g$.

Il r\'esulte du lemme \ref{calcul} que $\alpha_{C''}\in B^g$. 
Par cons\'equent, $\alpha_C\in B^g$, donc $\alpha_C\in B^{G'}$ puisque
$G'$ est engendr\'e par des doubles r\'eflexions
(corollaire \ref{parab2}). Donc, $\alpha_C\in B^G$ et le th\'eor\`eme
r\'esulte du lemme \ref{invar}.
\end{proof}

Le r\'esultat suivant pr\'ecise la premi\`ere implication du
th\'eor\`eme \ref{constr}.

\begin{thm}\label{theo 326}
Soit $G$ un bon sous-groupe de $\tW$ ne contenant pas de
r\'eflexions. Pour $H\in\CA(\tW)$, soit
$s_H$ le g\'en\'erateur de $\tW_H$ de d\'eterminant $\exp(2i\pi/|\tW_H|)$.

Alors, $G$ est engendr\'e par les
$s_{H'} s_H^{-1}$ o\`u $\{H,H'\}$ d\'ecrit les paires
d'hyperplans de $\CA'(\tW)$ ayant la m\^eme image dans
$\overline{\CA'}(\tW)$.
\end{thm}

\begin{proof}
Soit $G'$ le sous-groupe de $G$ engendr\'e par les $s_{H'}s_H^{-1}$.
Puisque $G\tW_H=G\tW_{H'}$, $\det(s_H)=\det(s_{H'})$ et
$G<SL(\tV)$, on a $G'\le G$. 

Si $G\tW_H=G\tW_{H'}$ pour $H,H'\in\CA(\tW)$, alors
$G'\tW_H=G'\tW_{H'}$. Par cons\'equent, il r\'esulte du lemme
\ref{invar} que $G'$ est bon dans $\tW$ et du th\'eor\`eme \ref{prop} (iv)
que $B^{G'}=B^G$. Par cons\'equent, $G'=G$.
\end{proof}

\medskip

\begin{rem}
Il ne suffit pas que $G$ soit engendr\'e par des \'el\'ements de la forme $s's^{-1}$ 
comme dans le th\'eor\`eme \ref{theo 326} pour qu'il soit bon, 
comme le montre l'exemple suivant.

Soit $\tW=B_2(4)$ (sous-groupe de $GL_2(\BC)$ des matrices monomiales
\`a coefficients non nuls dans $\mu_4(\BC)$) et $G$ le sous-groupe
engendr\'e par $g=-1$. Soit
$s=\left(\begin{array}{cc}
         0 & 1 \\
         1 & 0 
         \end{array}\right)$ 
et
$s'=\left(\begin{array}{cc}
          0 & -1 \\
         -1 &  0 
         \end{array}\right)$.
Alors, $g=s's^{-1}$ et $s$, $s'$ engendrent les groupes $\tW_{\ker s-1}$ et
$\tW_{\ker s'-1}$ respectivement. Cependant, $G$ n'est pas bon dans $\tW$ 
car il n'est pas bon dans le sous groupe de $\tW$ form\'e des matrices diagonales 
qui est \'egal \`a $\mu_4(\BC) \times \mu_4(\BC)$ (cf. propositions 
\ref{prop sous-groupe} et \ref{reduc}). 
\end{rem}


\bigskip

\subsection{Cas o\`u $W$ est ab\'elien} Nous \'etudions dans ce paragraphe le cas 
o\`u $W$ est ab\'elien, c'est-\`a-dire le cas o\`u $G$ contient le groupe des 
commutateurs de $\tW$. Dans ce cas, bien qu'un groupe ab\'elien puisse toujours \^etre 
r\'ealis\'e comme un groupe de r\'eflexion, il est possible que $W$ ne soit pas 
un groupe de r\'eflexion sur $V$, comme le montre l'exemple 
de la remarque \ref{cipasbon}. On a cependant les deux r\'esultats suivants~:

\begin{prop}\label{abelien}
Si $G=[\tW,\tW]$ est le sous-groupe engendr\'e par les commutateurs de $\tW$,
alors c'est un bon sous-groupe distingu\'e de $\tW$.
\end{prop}

\begin{proof}
Pour $H\in\CA(\tW)$, le groupe $G\tW_H$ est distingu\'e dans $\tW$. Par
cons\'equent, si $\CC\in\overline{\CA}(\tW)$, alors $\CC$ est $\tW$-stable,
donc $k\alpha_C$ est une repr\'esentation de dimension $1$ de $\tW$.
Alors, $G=[W,W]$ agit trivialement sur cette derni\`ere et
le lemme \ref{invar} montre que $G$ est bon dans $\tW$.
\end{proof}

Combin\'ee aux propositions \ref{quotient} et \ref{reduc}, la proposition
pr\'ec\'edente donne une description compl\`ete des cas o\`u $G$
est un bon sous-groupe distingu\'e de $\tW$ contenant $[\tW,\tW]$.

\begin{rem}
Stanley a d\'etermin\'e les cas o\`u $G=\tW\cap SL(\tV)$ est
d'intersection compl\`ete \cite[Th\'eor\`eme 5.1]{Sta}. Comme la
remarque \ref{cipasbon} le montre, $\tW\cap SL(\tV)$ peut \^etre
d'intersection compl\`ete sans \^etre un bon sous-groupe distingu\'e.
\end{rem}

\begin{prop}\label{pratique}
Soit $G=\tW \cap SL(\tV)$. On suppose que $\tW/G$ est d'ordre premier 
$p$. Alors $\overline{\CA'}(\tW)$ 
est r\'eduit \`a un seul \'el\'ement $C$. De plus, le groupe 
$G$ est bon dans $\tW$ et l'alg\`ebre $B^G$ est engendr\'ee 
par les invariants fondamentaux de $\tW$ et $\alpha_C$~: $\tV/G$ 
est donc une hypersurface dans $V$. L'id\'eal $I$ est engendr\'e 
par un seul \'el\'ement qui est obtenu en remarquant que $\alpha_C^p 
\in B^\tW$. 

Soient $d_1$,\dots, $d_n$ les degr\'es de $\tW$ et soit $N=|\CA(\tW)|$. 
Alors les degr\'es de $G$ sont $d_1$,\dots, $d_n$, $N$ et le degr\'e 
de la relation de $G$ est $Np$.
\end{prop}

\begin{proof} Cela r\'esulte du lemme \ref{invar} et de 
\cite[Th\'eor\`eme 3.1 et Corollaire 5.6]{Sta}.
\end{proof}

\bigskip

\section{Classification}~

\medskip

Le but de cette partie est de d\'eterminer les paires $(G,\tW)$ o\`u
$\tW$ est un groupe de r\'eflexion et $G$ un bon sous-groupe
distingu\'e de $\tW$. On supposera pour cela que $k=\BC$. 

\medskip
Le th\'eor\`eme \ref{constr} montre qu'il suffit de traiter le cas o\`u
$\tV$ est de dimension $2$.
Gr\^ace aux propositions \ref{quotient} 
et \ref{Grefl}, on peut remplacer $G$ et $\tW$ par $G/G_r$ et $\tW/G_r$,
o\`u $G_r$ est le sous-groupe de $G$ engendr\'e par les r\'eflexions de $G$~:
on \'etudiera donc le cas o\`u $G$ ne contient pas de r\'eflexion.  
D'apr\`es le corollaire \ref{propstrict}, cela implique que $G$ est 
contenu dans $SL(\tV)$.

Dans cette partie, nous ferons donc l'hypoth\`ese suivante~:

\medskip

{\em On suppose que $G$ est contenu dans $SL(\tV)$.
En outre, et ce jusqu'\`a la section \ref{classificationrangsuperieur},
on suppose que $\dim\tV=2$.}

\medskip

Une r\'ef\'erence possible pour les r\'esultats utilis\'es concernant
les groupes de r\'eflexion complexes est \cite{BrMaRou}.

\bigskip
Fixons tout d'abord quelques notations. Soit $\tV$ un 
espace vectoriel de dimension $2$ muni de sa graduation standard 
et soit $(x,y)$ une base de $\tV$. On note 
$(X,Y)$ la base de $\tV^*$ duale de $(x,y)$. Via le choix de cette base, 
nous identifierons $GL(\tV)$ avec $GL_2(\BC)$.

\bigskip

Il est bien connu que tout sous-groupe $G$ de $SL_2(\BC)$ est d'intersection 
compl\`ete et que son alg\`ebre d'invariants $B^G$ est engendr\'ee par trois 
polyn\^omes homog\`enes $p_1$, $p_2$ et $p_3$. On pose $d_i=\deg p_i$.
Les $p_i$ seront num\'erot\'es de sorte que 
$d_1 \le d_2 \le d_3$ et choisis de sorte que l'espace 
qu'ils engendrent dans $B^G$ soit stable sous l'action du normalisateur $N(G)$ 
de $G$ dans $GL_2(\BC)$. On notera $V^*$ cet espace vectoriel et on 
notera $(X_1,X_2,X_3)$ les \'el\'ements de $A=k[V]$ correspondant \`a
$(p_1,p_2,p_3)$~: l'image de $X_i$ par le morphisme surjectif
d'alg\`ebres gradu\'ees $A=k[X_1,X_2,X_3] \to B^G$ est $p_i$.
On notera $(x_1,x_2,x_3)$ 
la base de $V$ duale de $(X_1,X_2,X_3)$. Le groupe $GL(V)$ sera identifi\'e 
\`a $GL_3(\BC)$ par le choix de cette base. On notera 
$$\ph : N(G) \to GL_3(\BC)$$
le morphisme de groupes induit.

Puisque $G$ est d'intersection compl\`ete, le noyau $I$ du morphisme
$A \to B^G$ est engendr\'e par un \'el\'ement 
homog\`ene $R$ de degr\'e $e$ (le degr\'e de $X_i$ est $d_i$).
Avec ces notations, on a
$$|G|={d_1d_2d_3 \over e}.$$

\bigskip

\subsection{$G$ d'ordre $2$\label{sub C2}} 
\label{Gordre2}
Avant de s'int\'eresser au cas 
g\'en\'eral, nous allons d\'ecrire la situation lorsque 
$G$ est d'ordre $2$. Supposons donc que $G=\mu_2(\BC)$,
l'unique sous-groupe d'ordre $2$ de $SL_2(\BC)$ (on identifie
$\BC^\times$ avec $Z(GL_2(\BC))$).
On prend
$p_1=XY$, $p_2=X^2$, $p_3=Y^2$ et $R=X_1^2-X_2X_3$. On a 
donc $d_1=d_2=d_3=2$ et $e=4$ car tous les $X_i$ sont de degr\'es $2$. 
L'ensemble des degr\'es de $G$ est alors $\{2,2,2\}$ et l'ensemble des degr\'es 
des relations est $\{4\}$.

Le normalisateur de $G$ dans $GL_2(\BC)$ est $GL_2(\BC)$ et, si $g \in GL_2(\BC)$, 
alors $\ph(g)(R)=(\det g)^2 R$. Si $CO(R)$ est le groupe orthogonal conforme 
de la forme quadratique $R$ sur $V$, alors l'image du morphisme naturel 
de groupes alg\'ebriques $\ph : GL_2(\BC) \to GL(V)$ est contenu dans $CO(R)$. 
De plus,
$$N(G,rel)=\sqrt{SL}_2(\BC)=\{g \in GL_2(\BC)~|~g^2 \in SL_2(\BC)\}.$$
et l'image de $N(G,rel)$ par le morphisme $\ph$ est le groupe 
orthogonal $O(R)$ de la forme quadratique $R$. 
La proposition suivante r\'esulte de ce qui pr\'ec\`ede, du th\'eor\`eme \ref{ci} et 
de la proposition \ref{degres}~:

\begin{prop}\label{degres rang 2}
Soit $\tW$ un groupe de r\'eflexion complexe de rang $2$ et contenant $G=\{1,-1\}$. 
Alors $G$ est bon dans $\tW$ si et seulement si $\tW$ est engendr\'e 
par des r\'eflexions d'ordre $2$. Dans ce cas, si $2$, $a$ et $b$ sont les
degr\'es 
de $W$ dans $V$, alors les degr\'es de $\tW$ dans $\tV$ sont $2a$ et $2b$.
\end{prop}

Notons les r\'esultats \'el\'ementaires suivants~:

\begin{lemme}
Soit $g$ un \'el\'ement de $GL_2(\BC)$ dont les valeurs propres sont $\l$ et $\m$. 
Alors les valeurs propres de $\ph(g)$ dans $V$ sont $\l^2$, $\m^2$ et $\l\m$.
\end{lemme}

\begin{cor}
Si $\tilde{s}$ est une r\'eflexion d'ordre $2$ de $GL_2(\BC)$, alors $\ph(\tilde{s})$ 
est une r\'eflexion d'ordre $2$ de $O(R)$. R\'eciproquement, si $s$ est une 
r\'eflexion de $O(R)$, alors $\ph^{-1}(s)$ est form\'e de deux r\'eflexions 
d'ordre $2$.
\end{cor}

Le morphisme $\ph$ \'etablit une bijection entre l'ensemble des classes de
conjugaison de sous-groupes de 
$GL_2(\BC)$ engendr\'es par des r\'eflexions d'ordre $2$ et contenant $-1$ et 
l'ensemble des classes de conjugaison de sous-groupes non triviaux de $O(R)$
engendr\'es par des r\'eflexions.
 
Remarquons qu'un sous-groupe fini de $O(R)$ stabilise une structure 
r\'eelle sur l'espace vectoriel $V$, donc la correspondance ci-dessus d\'efinit
une bijection entre les classes de conjugaison de sous-groupes finis 
de $GL_2(\BC)$ engendr\'es par des r\'eflexions d'ordre $2$ et contenant 
$-1$ et l'ensemble des classes de conjugaison de sous-groupes finis non
triviaux de $GL_3(\BR)$ engendr\'es par des r\'eflexions, c'est-\`a-dire
l'ensemble des classes d'isomorphisme de groupes de Coxeter
de rang compris entre $1$ et $3$. 

Un cas particuli\`erement int\'eressant est celui o\`u le groupe $\tW/G$ est irr\'eductible. 
Il n'y a que trois classes d'isomorphisme de groupes de Coxeter irr\'eductibles de 
rang $3$, les groupes de type $A_3$, $B_3$ et $H_3$.
Leurs images r\'eciproques par $\ph$ sont les groupes de type $G_{12}$,
$G_{13}$ et $G_{22}$ respectivement.

La table $1$ donne la liste des groupes de r\'eflexion complexes
de rang $2$, 
engendr\'es par des r\'eflexions d'ordre $2$ et contenant $G=\mu_2(\BC)$,
leurs degr\'es, 
le type du groupe $\tW/G$ et les degr\'es de $\tW/G$~; $d$
d\'esignera un entier naturel non nul (on convient que $I_2(1)=A_1$).

\begin{table}
$\begin{array}{|c|c|c|c|}
\hline
\tW & {\mathrm{degr\acute{e}s~de~}}\tW & \tW/\{1,-1\} & 
{\mathrm{degr\acute{e}s~de~}}\tW/\{1,-1\} \vphantom{\DS{A \over B}}\\
\hline
G_{12} & 6,~8 & A_3 & 2,~3,~4 \vphantom{\DS{A \over B}}\\
G_{13} & 8,~12 & B_3 & 2,~4,~6 \vphantom{\DS{A \over B}}\\
G_{22} & 12,~20 & H_3 & 2,~6,~10 \vphantom{\DS{A \over B}}\\
I_2(2d) & 2,~2d & I_2(d) & 1,~2,~d
\vphantom{\DS{A \over B}} \\
~G(2d,d,2)~~ & 2d,~4 & I_2(d) \times A_1 & 2,~2,~d \vphantom{\DS{A \over B}}\\
\hline\end{array}$

\smallskip
\caption{}
\end{table}

\subsection{Sous-groupes de $SL_2(\BC)$ et $GL_2(\BC)$}
Si $d$ est un entier naturel non nul, on notera $\m_d$ le 
groupe $\m_d(\BC)$.
Pour $\z \in \BC^\times$, on pose~:
$$t(\z)=\matrice{\z & 0 \\ 0 & \z^{-1}}.$$
On pose aussi 
$$\s=\matrice{0 & -1 \\ 1 & 0}$$
$$s=\matrice{0 & 1 \\ 1 & 0}.\leqno{\mathrm{et}}$$

On pose 
$$C_d=\{t(\z)~|~\z^d=1\},$$
$$I_2(d)=< C_d,s>$$
$$\tilde{I}_2(d)=<\s,C_{2d} >.\leqno{\mathrm{et}}$$
Le groupe $\tilde{I}_2(d)$ est une extension centrale du groupe
di\'edral d'ordre $2d$, $I_2(d)$, par $\BZ/2\BZ$.
Notons que $\tilde{I}_2(1)$ est conjugu\'e \`a $C_4$ 
dans $GL_2(\BC)$.

On note $\tilde{\Sn}_4$ le normalisateur de $\tilde{I}_2(2)$ dans $SL_2(\BC)$. 
Soit $\AGt_4$ le groupe d\'eriv\'e de $\tilde{\Sn}_4$.
Le groupe $\SGt_4$ (respectivement $\AGt_4$) est une extension centrale du 
groupe sym\'etrique de degr\'e $4$ (respectivement du groupe altern\'e de 
degr\'e $4$) par $\BZ/2\BZ$.
On note $\AGt_5$ un sous-groupe de $SL_2(\BC)$ extension centrale du
groupe altern\'e de degr\'e $5$ par $\BZ/2\BZ$.
La proposition suivante est bien connue~:

\begin{prop}
Soit $G$ un sous-groupe fini de $SL_2(\BC)$. Alors $G$ est conjugu\'e, dans 
$GL_2(\BC)$, \`a un et un seul des groupes $C_d$ $(d \ge 1)$, $\tI_2(d)$ $(d \ge 2)$, 
$\AGt_4$, $\SGt_4$ et $\AGt_5$.
\end{prop}

Pour finir, on pose
$$T=\left\{ \matrice{a & 0 \\ 0 & b}~\Bigg\vert ~a,b \in \BC^\times\right\} ,$$
$$C_\infty=\{t(\z)~|~\z \in \BC^\times\}$$
$$I_2(\infty)=\langle C_\infty, s\rangle.\leqno{\mathrm{et}}$$
La table $2$
donne, en fonction du groupe $G$, les valeurs des param\`etres $|G|$, 
$d_1$, $d_2$, $d_3$, $e$, $N(G)$ et $N(G,rel)$.

\begin{table}
\begin{tabular}{|c|c|c|c|c|c|}
\hline
$G$ \vphantom{$\DS{A \over B}$} & $|G|$ & $d_1$, $d_2$, $d_3$ & $e$ & $N(G)$ & $N(G,rel)$ \\
\hline
$C_2$\vphantom{$\DS{A \over B}$} & $2$ & $2$, $2$, $2$ & $4$ & $GL(\tV)$ & $\sqrt{SL}(\tV)$ \\
$C_d$\vphantom{$\DS{A \over B}$}
 $(d \ge 3)$ & $d$ & $2$, $d$, $d$ & $2d$ & $\langle s,T \rangle$ &
 $I_2(\infty) \m_{2d}$ \\
$\tI_2(d)$\vphantom{$\DS{A \over B}$}
 $(d \ge 3)$ & $4d$ & $4$, $2d$, $2(d+1)$ & $4(d+1)$ & $\tI_2(2d)\BC^\times$ &
   $\tI_2(2d)\m_{4(d+1)}$ \\
$\tI_2(2)$ \vphantom{$\DS{A \over B}$}
& $8$ & $4$, $4$, $6$ & $12$ & $\SGt_4\BC^\times$ & $\SGt_4\m_{12}$ \\
$\AGt_4$ \vphantom{$\DS{A \over B}$}
& $24$ & $6$, $8$, $12$ & $24$ & $\SGt_4 \BC^\times$ & $\SGt_4\m_{24}$ \\
$\SGt_4$ \vphantom{$\DS{A \over B}$}
& $48$ & $8$, $12$, $18$ & $36$ & $\SGt_4\BC^\times$ & $\SGt_4\m_{36}$ \\
$\AGt_5$\vphantom{$\DS{A \over B}$}
 & $120$ & $12$, $20$, $30$ & $60$ & $\AGt_5\BC^\times$ & $\AGt_5 \m_{60}$ \\
\hline
\end{tabular}
\smallskip
\caption{}
\end{table}

\begin{rem} Le groupe $\tI_2(2)$ est le groupe des quaternions.
\end{rem}

Justifions maintenant le tableau pour $G=C_d$ et
$G=\tI_2(d)$ (les trois derniers cas se traitent
de mani\`ere similaire et sont laiss\'es au lecteur).

\medskip
\subsubsection{$G=C_d$\label{subsub Cd}}
Soit $G=C_d$ pour $d\ge 3$.
Le normalisateur $N(G)$ de 
$G$ dans $GL_2(\BC)$ est
$$N(G)=\langle T,s\rangle.$$

On prend $p_1=XY$, $p_2=X^d$, $p_3=Y^d$ et $R=X_1^d -X_2X_3$. 
On a donc $d_1=2$, $d_2=d_3=d$ et $e=2d$. On a $s\in N(G,rel)$, donc
$N(G,rel)=\langle s, N(G,rel)\cap T\rangle$. Finalement,
$$N(G,rel)=I_2(\infty)\m_d.$$

\medskip

\subsection{$G=\tI_2(d)$}
Consid\'erons maintenant $G=\tI_2(d)$.
On prend $p_1=X^2Y^2$, $p_2=X^{2d}+Y^{2d}$, $p_3=XY(X^{2d}-Y^{2d})$ 
et $R=X_3^2-X_1(X_2^2-4 X_1^d)$. On a donc $d_1=4$, $d_2=2d$, $d_3=2(d+1)$ et 
$e=4(d+1)$. Calculons maintenant le normalisateur de $G$. Tout d'abord, 
remarquons que le groupe d\'eriv\'e de $G$ est le groupe $C_d$ donc le 
normalisateur de $G$ est contenu dans le normalisateur de $C_d$. 

\medskip

Si $d \ge 3$, alors le normalisateur de $\tI_2(d)$ est contenu dans 
$<s,T>$.
Pour d\'eterminer $N(G)$ et $N(G,rel)$, il suffit de d\'eterminer 
leur intersection avec $T$, puisque $s\in N(G,rel)$ ~:
$$T\cap N(G)=\left\{ \matrice{a&0\\0&b} \Bigg\vert \ a^{2d}=b^{2d}\right\} $$
$$T\cap N(G,rel)=
 \left\{ \matrice{a & 0 \\ 0 &b} \Bigg\vert \ a^{2d}=b^{2d},\ (ab)^{2d+2}=1
 \right\}.$$

\medskip
Supposons maintenant $d=2$.
On a $N(G)=\SGt_4\cdot\BC^\times$ et, pour montrer 
que $N(G,rel)=\SGt_4\m_{12}$, il suffit de montrer que $\SGt_4$ est contenu
dans $N(G,rel)$ et d'utiliser la remarque \ref{norm}. Montrons donc que
$\SGt_4$ stabilise $R$. L'espace vectoriel des invariants de $G$ de degr\'e $6$ 
est de dimension $1$ et engendr\'e par $p_3$ et on note
$\e : \SGt_4 \to \BC^\times$ 
le caract\`ere lin\'eaire de $\SGt_4$ par lequel $\SGt_4$ agit sur cet 
espace vectoriel. Alors, $\SGt_4$ agit sur 
$\BC R$ via le caract\`ere lin\'eaire $\e^2$ de 
$\SGt_4$. Comme l'ab\'elianis\'e de $\SGt_4$ est d'ordre $2$, le caract\`ere
$\e^2$ est trivial, donc $\SGt_4 \subset N(G,rel)$.

\medskip

\begin{rem}
Si $G=\AGt_4$, $\SGt_4$ ou $\AGt_5$, alors les degr\'es de $G$ sont tous
distincts. 
Il en r\'esulte que si $\tW$ est un groupe de r\'eflexion de $GL_2(\BC)$ 
contenant $G$ comme bon sous-groupe distingu\'e, alors $\tW/G$ est ab\'elien. 
\end{rem}

\medskip

\subsection{Classification en rang deux.}
La table $3$ donne la liste, \`a conjugaison pr\`es, 
des groupes de r\'eflexion $\tW$ irr\'eductibles de rang $2$,
le type de $\tW\cap SL_2(\BC)$, de $\tW/Z(\tW)$ ainsi que
$a_\tW=[\tW:\tW\cap SL_2(\BC)]$.

\begin{table}
{\small
$\begin{array}{|c|c|c|c|}
\hline
\phantom{|^{|^|}_|} \tW \phantom{|^{|^|}_|} & \tW\cap SL_2(\BC) & \tW/Z(\tW) & a_\tW\\
\hline
G(mn,n,2)~  \vphantom{\DS{A \over B}} & \tI_2(mn/2) & &\\
(m\text{ pair}) & & &\\
G(mn,n,2)~  \vphantom{\DS{A \over B}} &C_{mn} & &\\
(m\text{ impair}) & & &\\
G_4 & \tI_2(2) & \AG_4 & 3\\
G_5 & \AGt_4 & \AG_4 & 3\\
G_6 & \tI_2(2) & \AG_4 & 6\\
G_7 & \AGt_4 & \AG_4 & 6\\
G_8 & \AGt_4 & \SG_4 & 4\\
G_9 & \SGt_4 & \SG_4 & 4\\
G_{10} & \AGt_4 & \SG_4 & 12\\
\hline\end{array}\quad\quad\quad
\begin{array}{|c|c|c|c|}
\hline
\phantom{|^{|^|}_|} \tW \phantom{|^{|^|}_|} & \tW\cap SL_2(\BC) & \tW/Z(\tW)& a_\tW\\
\hline
\qquad G_{11} \qquad & \SGt_4 & \SG_4 & 12\\
G_{12} & \AGt_4 & \SG_4 &2\\
G_{13} & \SGt_4 & \SG_4& 2\\
G_{14} & \AGt_4 & \SG_4 & 6\\
G_{15} & \SGt_4& \SG_4 & 6\\
G_{16} & \AGt_5 & \AG_5 & 5\\
G_{17} & \AGt_5 & \AG_5 & 10\\
G_{18} & \AGt_5 & \AG_5 & 15\\
G_{19} & \AGt_5 & \AG_5 & 30\\
G_{20} & \AGt_5 &\AG_5 & 3\\
G_{21} & \AGt_5 & \AG_5 &6\\
G_{22} & \AGt_5 & \AG_5& 2\\
\hline
\end{array}$
\smallskip
}
\caption{}
\end{table}

Pour la classification en dimension sup\'erieure, nous aurons aussi besoin
des groupes de r\'eflexions de rang $2$ non irr\'eductibles. Ceux-ci sont 
produits directs de deux groupes cycliques d'ordre $p$ et $q$.
Pour \'eviter la confusion
avec $C_d$, on pr\'ef\`erera la notation
$$\BZ/p\BZ \oplus \BZ/q\BZ :=\left\{ \begin{pmatrix}
\zeta_1 & 0 \\ 0 & \zeta_2 \end{pmatrix} \Bigg\vert
\zeta_1^p=\zeta_2^q=1 \right\}.$$
Soit $p\wedge q$ le plus grand diviseur commun de $p$ et $q$.
L'intersection de $\BZ/p\BZ \oplus \BZ/q\BZ$ avec $\SL_2(\BC)$
est $C_{p\wedge q}$.
Soit $d$ un diviseur de $p\wedge q$. Pour que $C_d$ soit bon dans 
$\BZ/p\BZ \oplus \BZ/q\BZ$, il faut et il suffit que
 $\BZ/p\BZ \oplus \BZ/q\BZ \subset N(C_d, rel) = I_2(\infty)\mu_{2d}$,
ce qui ne se produit que si $p|d$ et $q|d$. Ainsi les seules paires
$(\tW,G)$ o\`u
$\tW$ est r\'eductible et $G$ est non trivial et bon sont de la forme
$(\BZ/p\BZ \oplus \BZ/p\BZ,C_p)$. 

\bigskip

\bigskip

Les calculs effectu\'es pr\'ec\'edemment permettent maintenant d'obtenir
la classification~:

\begin{prop}
Soit $\tW$ un groupe de r\'eflexion de rang $2$ et soit $G$ un 
sous-groupe distingu\'e de $\tW$ contenu dans $\SL_2(\BC)$. Alors $G$ est bon dans $\tW$ 
si et seulement si on est dans l'un des cas suivants~:

\tete{a} $G=C_2$ et $\tW$ est engendr\'e par des r\'eflexions d'ordre $2$
(cf Table 1).

\tete{b} $G=C_d$ ($d \ge 3$) et $\tW=G(mn,n,2)$ avec $m|d$ et
$d|mn$.

\tete{c} $G=\tI_2(d)$ ($d \ge 3$) et $\tW=G(2d,d,2)$ ou $\tW=G(4d,4d,2)$.

\tete{d} $G=\tI_2(2)$ et $\tW=G(4,2,2)$ ou
                 $\tW=G_i$ avec $4 \le i \le 7$ ou $12 \le i \le 15$.

\tete{e} $G=\AGt_4$ et $\tW=G_i$, $i=5$ ou $7 \le i \le 15$.

\tete{f} $G=\SGt_4$ et $\tW=G_{15}$.

\tete{g} $G=\AGt_5$ et $\tW=G_i$, $16 \le i \le 22$.
\end{prop}

Pour $\tW=G(mn,n,2)$, on a utilis\'e le fait suivant~:
on a $G(mn,n,2)=<s,G(mn,n,2)\cap T>$ et 
$G(mn,n,2)\cap T=\{\matrice{a&0\\0&b}|\ a^{mn}=b^{mn}=1,\ (ab)^{m}=1\}$.

\smallskip
\'Etudions maintenant deux exemples de groupes $\tW$ exceptionnels.

Consid\'erons d'abord $\tW=G_{15}$. Puisque
$\tW\cap SL_2(\BC)=\SGt_4$, les seuls $G$ distingu\'es dans
$\tW$ sont $\mu_2(\BC)$ (qui a d\'ej\`a \'et\'e \'etudi\'e), $\tI_2(2)$,
$\AGt_4$ et $\SGt_4$ (ce sont tous des sous-groupes caract\'eristiques
de $\SGt_4$).

On a $\tW/Z(\tW)\simeq \SG_4$, donc $\tW\le \SGt_4\cdot\mu_\infty$.
Puisque $[\tW:\tW\cap SL_2(\BC)]=6$, on a plus pr\'ecis\'ement
$\tW\le \SGt_4\cdot \mu_{12}$. On en d\'eduit que  $\tI_2(2)$,
$\AGt_4$ et $\SGt_4$ sont bons dans $\tW$.

\smallskip
Pour $\tW=G_{11}$, on a aussi $\tW\cap SL_2(\BC)=\SGt_4$ et
$\tW\le \SGt_4\cdot\mu_\infty$. On a $[\tW:\tW\cap SL_2(\BC)]=12$, donc
$\tW\le \SGt_4\cdot \mu_{24}$. Puisque $|\tW|=|\SGt_4\cdot\mu_{24}|$, on
a l'\'egalit\'e $\tW=\SGt_4\cdot\mu_{24}$. On en d\'eduit que
$\AGt_4$ est bon dans $\tW$ mais que $\tI_2(2)$ et $\SGt_4$ ne sont
pas bons dans $\tW$.

\medskip

La table $4$ contient toutes les paires $(G,\tW)$ telles que $G$ est bon 
dans $\tW$, $G \not=\{1,-1\}$ et $W=\tW/G$ n'est pas ab\'elien. 

\begin{table}
$\begin{array}{|c|c|c|c|c|}
\hline
(G,\tW) \vphantom{\DS{A \over B}}
& {\mathrm{degr\acute{e}s~de~}}\tW & e & W & {\mathrm{degr\acute{e}s~de~}} W \\
\hline
(C_d,G(mn,n,2)) & mn,2m & 2d & G(m,1,1)\times G(mn/d,mn/d,2) & m,mn/d,2 \\
(\tI_2(2),G_{12}) \vphantom{\DS{A \over B}}& 6,~8 & ~~12~~ &A_2 & 2,~3,~1 \\
(\tI_2(2),G_{13}) \vphantom{\DS{A \over B}}& 8,~12 & ~~12~~ &A_2\times A_1 & 2,~3,~2 \\
(\tI_2(2),G_{14}) \vphantom{\DS{A \over B}}& 6,~24 & 12 &G(3,1,2) & 3,~6,~1 \\
(\tI_2(2),G_{15}) \vphantom{\DS{A \over B}}& 12,~24 & 12 &G(3,1,2) \times A_1 & 3,~6,~2 \\
\hline
\end{array}$
\smallskip
\caption{}
\end{table}

\subsection{Classification en rang sup\'erieur}
\label{classificationrangsuperieur}

Gr\^ace au th\'eor\`eme \ref{constr}, le probl\`eme g\'en\'eral de la
classification des paires $(\tW,G)$ se ram\`ene \`a la classification
en rang $2$ \'etablie ci-dessus.
Pr\'ecisons comment cette r\'eduction s'op\`ere. Soit
$\tW\subset\GL(V)$ un groupe
de r\'eflexion et soit $G$ bon dans $\tW$.
D'apr\`es \ref{constr} (iii), si $G$ est non trivial, il existe un sous-espace
$L$ de $V$ de codimension $2$ et un sous-groupe non trivial $G'$ de $G$ tel que
$G'$ est un bon sous-groupe distingu\'e de $\tW_L$. La
classification en rang $2$ donne, pour chaque $\tW_L$, une liste
des $G'$ possibles. Soit $\overline{G'}$
la cl\^oture distingu\'ee de $G'$ dans $\tW$. Toujours gr\^ace
au th\'eor\`eme \ref{constr}, il est clair que $G'$ est bon dans
$\tW$. D'apr\`es 
\ref{quotient}, le groupe $G/\overline{G'}$
est bon dans $\tW/\overline{G'}$. Comme l'ordre de $G/\overline{G}'$
est strictement plus petit que celui de $G$, on peut, quitte \`a it\'erer 
l'op\'eration, reconstruire $G$ \`a partir de groupes
$G'$ de rang $2$.

On est ainsi ramen\'e \`a classifier, pour $\tW$ donn\'e, les paires
$(L,G')$ avec $L$ de codimension $2$ et $G'$ bon dans $\tW_L$, \`a conjugaison
pr\`es.
Il suffit d'ailleurs de consid\'erer, parmi ces paires, celles pour
lesquelles $G'$ est minimal, puisqu'il s'agit simplement d'initier le processus
r\'ecursif~; les tables que nous donnons plus bas en donnent la liste.
Les diagrammes ``\`a la Coxeter'' fournis dans \cite{BrMaRou} facilitent
la recherche des paires $(L,G')$. Sauf pour les groupes exceptionnels $G_{27}$,
$G_{29}$, $G_{33}$ et $G_{34}$,
ces diagrammes poss\`edent la propri\'et\'e suivante~: soit $S$
un syst\`eme de g\'en\'erateurs de $\tW$ correspondant au diagramme~;
pour tout
sous-groupe ``parabolique'', c'est-\`a-dire de
la forme $\tW_L$ pour un certain sous-espace $L$, il existe un
sous-ensemble $I\subset S$ tel que le sous-groupe engendr\'e par $I$
soit conjugu\'e \`a $\tW_L$.

\subsection{``Morphismes de diagrammes'' et nouveaux diagrammes}
Consid\'erons la pr\'esenta\-tion $<S | R>$ associ\'ee au diagramme standard
d'un groupe de r\'eflexion 
$\tW \subset \GL(V)$. Soit $I\subset S$. Le sous-groupe de $\tW$ engendr\'e par
les r\'eflexions de $I$ est de la forme $\tW_L$, o\`u $L$ est
un sous-espace de $V$. L'espace $L$ est
de codimension $2$ quand $I$ est compos\'e de $2$ \'el\'ements li\'es par
une relation de tresses usuelle, ou compos\'e de $3$ \'el\'ements li\'es
par une relation de valence $3$, par exemple une relation circulaire
$stu=tus=ust$ (voir \cite{BrMaRou}). Quotienter par un $G'$ bon dans
$\tW_L$ correspond \`a ajouter \`a $R$ une ou plusieurs 
relations suppl\'ementaires, liant les \'el\'ements de $I$. Soit $R'$ le
nouvel ensemble de relations. Le groupe quotient $W$ admet la pr\'esentation
$<S | R'>$.

Il se trouve que, dans la plupart des cas (en fait
dans tous les cas, si l'on consid\`ere les diagrammes modifi\'es que nous
proposons plus bas), cette pr\'esentation
est identique, pour peu que l'on \'elimine les g\'en\'erateurs et les 
relations redondants,
\`a celle donn\'ee par le diagramme standard de $W$.
Ainsi, le morphisme
$\tW \twoheadrightarrow W$ provient d'un ``morphisme de diagrammes'' entre
leurs diagrammes standard. Nous n'abordons pas cette formalisation dans
le pr\'esent travail. Il serait souhaitable de comprendre g\'eom\'etriquement
la cat\'egorie des groupes de r\'eflexions avec ``morphismes de
diagrammes''. Notons toutefois qu'elle exigerait une d\'efinition
pr\'ecise non seulement de la cat\'egorie des diagrammes, ce qui ne pose
pas de difficult\'e,
mais aussi du ``foncteur'' qui \`a un groupe de r\'eflexion associe son
diagramme (ou vice-versa). Deux difficult\'es apparaissent~: la pr\'esence
d'\'el\'ements de $N_{GL(\tV)}(\tW)/\tW$ qui ne peuvent pas \^etre
r\'ealis\'es comme automorphismes de diagrammes et la pr\'esence de
sous-groupes paraboliques qui ne sont pas conjugu\'es \`a des
sous-groupes paraboliques standard.

\smallskip
Les seules exceptions \`a la compatibilit\'e entre diagrammes sont
$G_{12}$, $G_{13}$ et $G_{22}$~: selon la table $1$, ils
admettent comme quotients de r\'eflexion les groupes de Coxeter
de type respectivement $A_3$, $B_3$ et $H_3$~; 
or, partant des diagrammes de \cite{BrMaRou}, on ne retrouve pas au
quotient des syst\`emes g\'en\'erateurs de Coxeter.
Nous proposons de nouveaux diagrammes, compatibles au quotient avec
les pr\'esentations de Coxeter.

Soient $a,b,c,e,f$ des entiers.  Le diagramme $\tI_{a,b,c}(e,f)$

$$
\xy (0,2) *++={a} *\frm{o} ; (13,2) *++={b} *\frm{o} **@{~} ;
 (26,2) *++={c} *\frm{o} **@{~}; (6,5) *++={e} ;
   (0,-3) *++={s} ; (13,-3) *++={t} ; (26,-3) *++={u}; (20,5) *++={f} \endxy
$$

symbolise la pr\'esentation
$$<s,t,u | s^a=t^b=u^c=1, 
\underbrace{\dots tstus}_{\text{$e+1$ termes}} =
\underbrace{\dots ststu}_{\text{$e+1$ termes}}, \; 
\underbrace{ustut \dots }_{\text{$f+1$ termes}} =
\underbrace{stutu \dots }_{\text{$f+1$ termes}} >.$$
La relation entre $s$ et $t$ est la relation de tresses de longueur
$e$ ``tordue'' en intercalant $u$ en derni\`ere et en avant-derni\`ere
position. De m\^eme pour celle entre $t$ et $u$, d'o\`u les ar\^etes
``tordues'' du diagramme.

Les groupes $G_{12}$, $G_{13}$ et $G_{22}$ correspondent \`a $a=b=c=2$ et
respectivement \`a $(e,f)$ \'egal \`a $(3,3)$, $(4,3)$ et $(5,3)$.
On v\'erifie ais\'ement qu'en supprimant les relations quadratiques on obtient
une pr\'esentation du groupe de tresses associ\'e.

Les groupes de r\'eflexion complexes de rang $2$, non r\'eels, contenant $-1$
et
engendr\'es par des r\'eflexions d'ordre $2$ correspondent aux paires
$(e,f)\in \{(3,3),(4,3),(5,3),(d,2) (d\ge 2)\}$.
Ce sont
exactement les paires $(e,f)$ avec $e,f\ge 2$ telles que $1/e+1/f>1/2$.

En ajoutant la relation $su=us$, les relations se simplifient et
l'on retrouve la pr\'esentation de Coxeter associ\'ee \`a
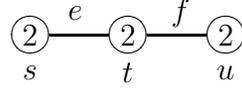
$$\xy (0,2) *++={2} *\frm{o} ; (13,2) *++={2} *\frm{o} **@{-} ;
 (26,2) *++={2} *\frm{o} **@{-}; (6,5) *++={e} ;
   (0,-3) *++={s} ; (13,-3) *++={t} ; (26,-3) *++={u}; (20,5) *++={f} \endxy
$$
(cf Table 1).

Le diagramme $\tI_{a,b,c}(e,f)$ redonne les pr\'esentations de \cite{BrMaRou}
pour un certain nombre de groupes de dimension $2$
(en dehors des cas $G_{12}$, $G_{13}$ et $G_{22}$ o\`u les pr\'esentations
sont nouvelles). Ils fournissent ainsi des diagrammes pour tous les groupes
de r\'eflexion complexes de rang $2$ qui ne sont pas de Shephard~:

\begin{itemize}
\item $I_{2,2,d}(e,2)$ pour $G(de,e,2)$, $d,e\ge 2$
\item $I_{2,2,3}(2,4)$ pour $G_{15}$
\item $I_{2,3,3}(2,2)$ pour $G_7$, 
$I_{2,3,4}(2,2)$ pour $G_{11}$ et $I_{2,3,5}(2,2)$ pour $G_{19}$.
\end{itemize}

\subsection{Tables des quotients \'el\'ementaires}
Soit $\tW \subset \GL(\tV)$ un groupe de
r\'eflexion complexe. Les tables $5$ et $6$ permettent
de retrouver tous les $G\subset \tW$ tels que $G\subset \SL(\tV)$ et
$G$ est bon dans $\tW$, ainsi que les quotients $W= \tW/G$ correspondants.
Consid\'erons le diagramme standard associ\'e \`a $\tW$, c'est-\`a-dire
la r\'eunion disjointe des diagrammes standard associ\'es aux facteurs
de r\'eflexion irr\'eductibles de $\tW$.
Les quotients $W= \tW/G$ sont les groupes de r\'eflexion obtenus
en appliquant une ou plusieurs des transformations \'el\'ementaires d\'ecrites.
(Les tables contiennent toutes les paires $(\tW,G)$ o\`u $G$ est bon, non
trivial, minimal --- pour simplifier les tables, nous n'avons pas exclus
certains cas o\`u $G$ n'est pas minimal).

Les tables se lisent de la fa\c con suivante. La premi\`ere colonne
donne le nom usuel d'un $\tW$ irr\'eductible, la seconde son diagramme.
La troisi\`eme colonne donne le diagramme quotient, obtenu en rajoutant
la relation donn\'ee dans la derni\`ere colonne. Un g\'en\'erateur
de $\tW$ est envoy\'e sur le g\'en\'erateur de m\^eme nom de $W$.
Comme il est expliqu\'e en \ref{classificationrangsuperieur}, un
quotient \'el\'ementaire est donn\'e par le choix d'un sous-espace
$L$ de $V$ de codimension $2$ et de $G'$ minimal bon dans $\tW_L$.

\begin{table}
{\small
$\begin{array}{|c|c|c|c|}
\hline
\phantom{\widetilde{|^|}}
\text{$\tW$} \phantom{\widetilde{|^|}} & \text{Diagramme de $\tW$} &
  \text{Diagramme de $W$}  & \text{Relation} \\
\hline
& & & \\
&\begin{array}{c} 
\text{Engendr\'e par des r\'eflexions } \\
s_1,\ldots, s_r \text{ d'ordre } p \end{array} & 
\xy  (10,2) *++={p} *\frm{o} ;
    (10,-3) *++={s_1} \endxy &
s_1=\cdots=s_r \\
\hline
& & & \\
\mathfrak{S}_4,G_{25} & 
\xy (0,2) *++={p} *\frm{o} ; (13,2) *++={p} *\frm{o} **@{-} ;
 (26,2) *++={p} *\frm{o} **@{-};
   (0,-3) *++={s} ; (13,-3) *++={t} ; (26,-3) *++={u} \endxy &
\xy (0,2) *++={p} *\frm{o} ; (13,2) *++={p} *\frm{o}
 **@{-};
   (0,-3) *++={s} ; (13,-3) *++={t} \endxy  &
   s=u \\
\hline
& & &  \\
G_{26} & 
\xy (0,2) *++={2} *\frm{o} ; (13,2) *++={3} *\frm{o} **@{=} ;
 (26,2) *++={3} *\frm{o} **@{-};
   (0,-3) *++={s} ; (13,-3) *++={t} ; (26,-3) *++={u} \endxy &
\xy (0,2) *++={2} *\frm{o} ; (13,2) *++={3} *\frm{o} ;
 (26,2) *++={3} *\frm{o} **@{-};
   (0,-3) *++={s} ; (13,-3) *++={t} ; (26,-3) *++={u} \endxy &
   st=ts \\
& & &  \\
\begin{array}{c} G_{28} \\ (=W(F_4)) \end{array} & 
\xy (0,2) *++={2} *\frm{o} ; (10,2) *++={2} *\frm{o} **@{-} ;
 (20,2) *++={2} *\frm{o} **@{=}; (30,2) *++={2} *\frm{o} **@{-} ;
   (0,-3) *++={s} ; (10,-3) *++={t} ; (20,-3) *++={u};
    (30,-3) *++={v} \endxy &
\xy (0,2) *++={2} *\frm{o} ; (10,2) *++={2} *\frm{o} **@{-} ;
 (20,2) *++={2} *\frm{o} ; (30,2) *++={2} *\frm{o} **@{-} ;
   (0,-3) *++={s} ; (10,-3) *++={t} ; (20,-3) *++={u};
    (30,-3) *++={v} \endxy &
 tu=ut 
\\
G_{29} &
\xy (0,2) *++={2} *\frm{o} ; (12,2) *++={2} *\frm{o} **@{-} ;
 (24,2) *++={2} *\frm{o} **@{=}; (18,9) *++={2} *\frm{o} **@{-};
 (12,2) *++={2} *\frm{o} ; (18,9) *++={2} *\frm{o} **@{-};
 (18,2) *++={}  ; (18,9) *++={2} *\frm{o} **@{=};
   (0,-3) *++={s} ; (12,-3) *++={t} ; (24,-3) *++={u};
    (23,9) *++={v} \endxy &
\xy (0,2) *++={2} *\frm{o} ; (10,2) *++={2} *\frm{o} **@{-} ;
 (20,2) *++={2} *\frm{o} **@{-}; (30,2) *++={2} *\frm{o} **@{-} ;
   (0,-3) *++={s} ; (10,-3) *++={t} ; (20,-3) *++={v};
    (30,-3) *++={u} \endxy &  tu=ut\\
& & &  \\
G_{31} &
\xy  
(2,7) *++={s} ; (25,7) *++={w} ; (6,7) *++={2} *\frm{o} ;
(0,-3) *++={2} *\frm{o} **@{-}; (13,-3) *++={2} *\frm{o} **@{-} ;
 (26,-3) *++={2} *\frm{o} **@{-}; (20,7) *++={2} *\frm{o} **@{-} ;
   (0,-8) *++={t} ; (13,-8) *++={u} ; (26,-8) *++={v} ;
   (13,4.3) *++={\phantom{sdds}} *\frm{o} \endxy
&
\xy  
(2,7) *++={s} ; (25,7) *++={w} ; (6,7) *++={2} *\frm{o} ;
(0,-3) *++={2} *\frm{o} **@{-}; (13,-3) *++={2} *\frm{o} **@{-} ;
 (26,-3) *++={2} *\frm{o} **@{-}; (20,7) *++={2} *\frm{o} **@{-} ;
   (0,-8) *++={t} ; (13,-8) *++={u} ; (26,-8) *++={v}
  \endxy
& su=us \\
\hline
& & &  \\
G(p,1,n) & 
\xy (0,2) *++={p} *\frm{o} ; (10,2) *++={2} *\frm{o} **@{=} ;
 (20,2) *++={2} *\frm{o} **@{-}; (30,2) *++={\dots} **@{-} ;
 (40,2) *++={2} *\frm{o} **@{-};
   (0,-3) *++={s} ; (10,-3) *++={t_1} ; (20,-3) *++={t_2};
    (40,-3) *++={t_{n-1}} \endxy &
\xy (0,2) *++={p} *\frm{o} ; (10,2) *++={2} *\frm{o} ;
 (20,2) *++={2} *\frm{o} **@{-}; (30,2) *++={\dots} **@{-} ;
 (40,2) *++={2} *\frm{o} **@{-};
   (0,-3) *++={s} ; (10,-3) *++={t_1} ; (20,-3) *++={t_2};
    (40,-3) *++={t_{n-1}} \endxy &
st_1=t_1s \\
G(md,md,n) &
\xy   (10,8) *++={2} *\frm{o}; (10,-4) *++={2} *\frm{o} **@{-};
     (9,2) *++={} ;  (20,2) *++={2} *\frm{o} **@{=};
   (10,8) *++={2} *\frm{o} ;  (20,2) *++={2} *\frm{o} **@{-};
  (10,-4) *++={2} *\frm{o} ;  (20,2) *++={2} *\frm{o} **@{-};
 (30,2) *++={\dots} **@{-} ; (40,2) *++={2} *\frm{o} **@{-};
   (5,8) *++={t_1} ; (5,-4) *++={t'_1} ; (20,-3) *++={t_2};
    (40,-3) *++={t_{n-1}}; (6,2) *++={md} \endxy 
    &
\xy   (10,8) *++={2} *\frm{o}; (10,-4) *++={2} *\frm{o} **@{-};
     (9,2) *++={} ;  (20,2) *++={2} *\frm{o} **@{=};
   (10,8) *++={2} *\frm{o} ;  (20,2) *++={2} *\frm{o} **@{-};
  (10,-4) *++={2} *\frm{o} ;  (20,2) *++={2} *\frm{o} **@{-};
 (30,2) *++={\dots} **@{-} ; (40,2) *++={2} *\frm{o} **@{-};
   (5,8) *++={t_1} ; (5,-4) *++={t'_1} ; (20,-3) *++={t_2};
    (40,-3) *++={t_{n-1}} ; (7,2) *++={m} \endxy 
& \begin{array}{c} \underbrace{t_1t'_1t_1\dots}_{\text{$m$ termes}} \\
 = \underbrace{t'_1t_1t'_1\dots}_{\text{$m$ termes}} \end{array} \\
G(d,d,n) &
\xy   (10,8) *++={2} *\frm{o}; (10,-4) *++={2} *\frm{o} **@{-};
     (9,2) *++={} ;  (20,2) *++={2} *\frm{o} **@{=};
   (10,8) *++={2} *\frm{o} ;  (20,2) *++={2} *\frm{o} **@{-};
  (10,-4) *++={2} *\frm{o} ;  (20,2) *++={2} *\frm{o} **@{-};
 (30,2) *++={\dots} **@{-} ; (40,2) *++={2} *\frm{o} **@{-};
   (5,8) *++={t_1} ; (5,-4) *++={t'_1} ; (20,-3) *++={t_2};
    (40,-3) *++={t_{n-1}}; (6,2) *++={d} \endxy 
    &
\xy  (10,2) *++={2} *\frm{o} ;
 (20,2) *++={2} *\frm{o} **@{-}; (30,2) *++={\dots} **@{-} ;
 (40,2) *++={2} *\frm{o} **@{-};
    (10,-3) *++={t_1} ; (20,-3) *++={t_2};
    (40,-3) *++={t_{n-1}} \endxy  & t_1=t'_1 \\
G(3d,3d,n) &
\xy   (10,8) *++={2} *\frm{o}; (10,-4) *++={2} *\frm{o} **@{-};
     (9,2) *++={} ;  (20,2) *++={2} *\frm{o} **@{=};
   (10,8) *++={2} *\frm{o} ;  (20,2) *++={2} *\frm{o} **@{-};
  (10,-4) *++={2} *\frm{o} ;  (20,2) *++={2} *\frm{o} **@{-};
 (30,2) *++={\dots} **@{-} ; (40,2) *++={2} *\frm{o} **@{-};
   (5,8) *++={t_1} ; (5,-4) *++={t'_1} ; (20,-3) *++={t_2};
    (40,-3) *++={t_{n-1}}; (6,2) *++={3d} \endxy 
    &
\xy  (10,2) *++={2} *\frm{o} ;
 (20,2) *++={2} *\frm{o} **@{-}; (30,2) *++={\dots} **@{-} ;
 (40,2) *++={2} *\frm{o} **@{-};
    (10,-3) *++={t'_1} ; (20,-3) *++={t_2};
    (40,-3) *++={t_{n-1}} \endxy  & t_1=t_2 \\
%
G(de,e,n) & 
\xy (20,10) *+\hbox{} ; (10,10) *+\hbox{}; **@{};
   (10,10) *+\hbox{}="p" ; (10,-6) *+\hbox{}="c"; {\ellipse_{-}};
   (-0.5,2) *++={d} *\frm{o} ; 
     (9,2) *++={} ;  (20,2) *++={2} *\frm{o} **@{=};
   (10,8) *++={2} *\frm{o} ;  (20,2) *++={2} *\frm{o} **@{-};
  (10,-4) *++={2} *\frm{o} ;  (20,2) *++={2} *\frm{o} **@{-};
 (30,2) *++={\dots} **@{-} ; (40,2) *++={2} *\frm{o} **@{-};
   (14,10) *++={t_1} ; (14,-6) *++={t'_1} ; (20,-3) *++={t_2};
    (40,-3) *++={t_{n-1}} ; (0,-5) *++={e+1}; (-5,2) *++={s} \endxy 
&
\xy   (0,2) *++={d} *\frm{o} ; 
   (10,8) *++={2} *\frm{o}; (10,-4) *++={2} *\frm{o} **@{-};
     (9,2) *++={} ;  (20,2) *++={2} *\frm{o} **@{=};
   (10,8) *++={2} *\frm{o} ;  (20,2) *++={2} *\frm{o} **@{-};
  (10,-4) *++={2} *\frm{o} ;  (20,2) *++={2} *\frm{o} **@{-};
 (30,2) *++={\dots} **@{-} ; (40,2) *++={2} *\frm{o} **@{-};
   (5,8) *++={t_1} ; (5,-4) *++={t'_1} ; (20,-3) *++={t_2};
    (40,-3) *++={t_{n-1}} ; (7,2) *++={e}; (0,-3) *++={s} \endxy 
& st_1=t_1s \\

\hline
\end{array}$
}
\smallskip

\caption{(rang sup\'erieur \`a $3$)}
\end{table}

{\bf \flushleft Exemple.} Le groupe $G(4,2,4)$ admet de multiples
quotients de r\'eflexion:

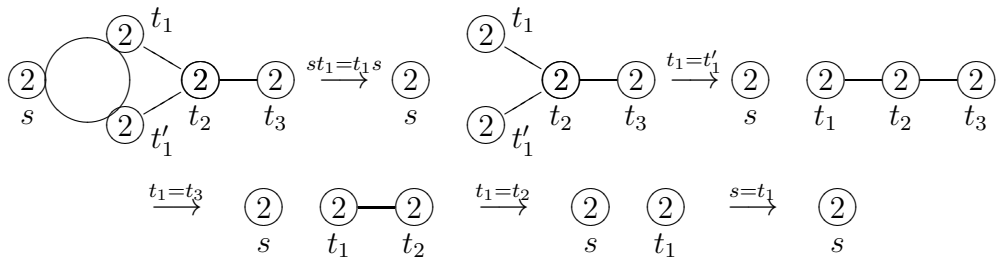
$$
\xy  (-3,0) *++={2} *\frm{o} ; (5,0) *++={\phantom{ww;.}} *\frm{o};
   (10,6) *++={2} *\frm{o} ;  (20,0) *++={2} *\frm{o} **@{-};
  (10,-6) *++={2} *\frm{o} ;  (20,0) *++={2} *\frm{o} **@{-};
 (30,0) *++={2} *\frm{o} **@{-};
   (-3,-5) *++={s}; (15,8) *++={t_1} ; (15,-8) *++={t'_1} ; (20,-5) *++={t_2};
    (30,-5) *++={t_3} \endxy \stackrel{st_1=t_1s}{\longrightarrow}
\xy  (0,0) *++={2} *\frm{o} ;
   (10,6) *++={2} *\frm{o} ;  (20,0) *++={2} *\frm{o} **@{-};
  (10,-6) *++={2} *\frm{o} ;  (20,0) *++={2} *\frm{o} **@{-};
 (30,0) *++={2} *\frm{o} **@{-};
   (0,-5) *++={s}; (15,8) *++={t_1} ; (15,-8) *++={t'_1} ; (20,-5) *++={t_2};
    (30,-5) *++={t_3} \endxy \stackrel{t_1=t'_1}{\longrightarrow}
\xy  (0,0) *++={2} *\frm{o} ;
   (10,0) *++={2} *\frm{o} ;  (20,0) *++={2} *\frm{o} **@{-};
 (30,0) *++={2} *\frm{o} **@{-};
   (0,-5) *++={s}; (10,-5) *++={t_1} ; (20,-5) *++={t_2};
    (30,-5) *++={t_3} \endxy 
    $$
$$    
\stackrel{t_1=t_3}{\longrightarrow} \quad \xy  (0,0) *++={2} *\frm{o} ;
   (10,0) *++={2} *\frm{o} ;  (20,0) *++={2} *\frm{o} **@{-};
   (0,-5) *++={s}; (10,-5) *++={t_1} ; (20,-5) *++={t_2};
   \endxy \quad \stackrel{t_1=t_2}{\longrightarrow} \quad
\xy  (0,0) *++={2} *\frm{o} ;
   (10,0) *++={2} *\frm{o} ;  (0,-5) *++={s}; (10,-5) *++={t_1} 
   \endxy  \quad \stackrel{s=t_1}{\longrightarrow} \quad
   \xy  (0,0) *++={2} *\frm{o} ;
   (0,-5) *++={s}
   \endxy
$$

\begin{table}
{\small
$\begin{array}{|c|c|c|c|}
\hline
\phantom{\widetilde{|^|}}
\text{$\tW$} \phantom{\widetilde{|^|}} & \text{Diagramme de $\tW$} &
  \text{Diagramme de $W$}  & \text{Relation} \\
\hline
& & & \\
I_{p,p}(m) &
\xy (0,2) *++={p} *\frm{o} ; (13,2) *++={p} *\frm{o} **@{-} ;
   (0,-3) *++={s} ; (13,-3) *++={t} ; (6,5) *++={m} \endxy &
\xy (0,2) *++={p} *\frm{o}; (0,-3) *++={s} \endxy &
 s=t \\

I_{p,q}(md) &
\xy (0,2) *++={p} *\frm{o} ; (13,2) *++={q} *\frm{o} **@{-} ;
   (0,-3) *++={s} ; (13,-3) *++={t} ; (6,5) *++={md} \endxy &
\xy (0,2) *++={p} *\frm{o} ; (13,2) *++={q} *\frm{o} **@{-};
   (0,-3) *++={s} ; (13,-3) *++={t} ; (6,5) *++={m} \endxy &
\begin{array}{c} \underbrace{sts\dots}_{\text{$m$ termes}} \\
 = \underbrace{tst\dots}_{\text{$m$ termes}} \end{array} \\

\tI_{a,b,c}(e,f) &
\xy (0,2) *++={a} *\frm{o} ; (13,2) *++={b} *\frm{o} **@{~} ;
 (26,2) *++={c} *\frm{o} **@{~}; (6,5) *++={e} ;
   (0,-3) *++={s} ; (13,-3) *++={t} ; (26,-3) *++={u}; (20,5) *++={f} \endxy

& 
\xy (0,2) *++={a} *\frm{o} ; (13,2) *++={b} *\frm{o} **@{-} ;
 (26,2) *++={c} *\frm{o} **@{-}; (6,5) *++={e} ;
   (0,-3) *++={s} ; (13,-3) *++={t} ; (26,-3) *++={u}; (20,5) *++={f} \endxy
& su=us \\
& & & \\
\hline
\end{array}$
}
\smallskip
\caption{(rang $2$)}
\end{table}

\begin{rem}
\begin{itemize}
\item[]
\item La transformation de $\mathfrak{S}_4$ (ou $G_{25}$) 
vers $\mathfrak{S}_3$ (ou $G_{4}$)
consiste \`a ``replier'' le diagramme. Une telle transformation pour
$\mathfrak{S}_n$ n'est pas possible si $n\ge 5$
(c'est li\'e \`a la simplicit\'e du groupe altern\'e $\mathfrak{A}_n$ pour
$n\ge 5$).
Ainsi tous les quotients non triviaux de $\mathfrak{S}_n$
s'obtiennent par des morphismes de diagrammes.
\item Excluons la r\`egle $\BZ/p\BZ \oplus \BZ/p\BZ \twoheadrightarrow
\BZ/p\BZ$. Partant de $\tW$ quelconque, et effectuant les transformations
au hasard, on aboutit toujours \`a un diagramme qui ne peut plus
\^etre transform\'e. Ce diagramme est celui de
$\tW^{\ab}= \tW/[\tW,\tW]$, qui est un produit direct de groupes de
r\'eflexion de rang $1$. 
\end{itemize}
\end{rem}

\subsection{Un exemple~: $G_{31}$}

Soient

$$s=\left( \begin{array}{cccc} 
 -1 & 0 & 0 & 0 \\
  0 & 1 & 0 & 0 \\
  0 & 0 & 1 & 0 \\
  0 & 0 & 0 & 1 
  \end{array} \right),\ \ 
t=\left( \begin{array}{cccc} 
  1/2 & -1/2 & -1/2 & -1/2 \\
 -1/2 &  1/2 & -1/2 & -1/2 \\
 -1/2 & -1/2 &  1/2 & -1/2 \\
 -1/2 & -1/2 & -1/2 &  1/2 
  \end{array} \right),\ \ 
u=\left( \begin{array}{cccc} 
 0 & i & 0 & 0 \\
-i & 0 & 0 & 0 \\
 0 & 0 & 1 & 0 \\
 0 & 0 & 0 & 1 
  \end{array} \right)$$
$$
v=\left( \begin{array}{cccc} 
 1 & 0 & 0 & 0 \\
 0 & 0 & 1 & 0 \\
 0 & 1 & 0 & 0 \\
 0 & 0 & 0 & 1 
  \end{array} \right),\ \
w=\left( \begin{array}{cccc} 
 0 & 1 & 0 & 0 \\
 1 & 0 & 0 & 0 \\
 0 & 0 & 1 & 0 \\
 0 & 0 & 0 & 1 
  \end{array} \right)$$

Soit $\tW$ le sous-groupe de $GL_4(\BC)$ engendr\'e par $s,t,u,v,w$.
C'est le groupe de r\'eflexion complexe $G_{31}$.

Soit $G=O_2(\tW)$, le plus grand $2$-sous-groupe distingu\'e de $\tW$.
C'est un groupe d'ordre $64$, engendr\'e par

$$\left( \begin{array}{cccc} 
 -1 & 0 & 0 & 0 \\
  0 & -1 & 0 & 0 \\
 0 & 0 & 1 & 0 \\
 0 & 0 & 0 & 1 
  \end{array} \right) ,
\left( \begin{array}{cccc}
 -1 & 0 & 0 & 0 \\
  0 & 1 & 0 & 0 \\
 0 & 0 & -1 & 0 \\
 0 & 0 & 0 & 1
  \end{array} \right) ,
\left( \begin{array}{cccc}
 0 & 1 & 0 & 0 \\
  1 & 0 & 0 & 0 \\
 0 & 0 & 0 & -1 \\
 0 & 0 & -1 & 0
  \end{array} \right) ,$$
$$\left( \begin{array}{cccc}
 0 & 0 & 1 & 0 \\
  0 & 0 & 0 & -1 \\
 1 & 0 & 0 & 0 \\
 0 & -1 & 0 & 0
  \end{array} \right) ,
\left( \begin{array}{cccc}
 0 & 0 & 0 & -i \\
  0 & 0 & -i & 0 \\
 0 & i & 0 & 0 \\
 i & 0 & 0 & 0
  \end{array} \right)$$ 

Les droites de la base canonique de $\BC^4$ 
forment un syst\`eme d'imprimitivit\'e de $G$~; le groupe sous-jacent des 
permutations de ces quatre droites est le sous-groupe kleinien $\frak K$ de 
$\mathfrak{S}_4$ (le sous-groupe engendr\'e par les doubles transpositions). 
Les coefficients non nuls de la matrice d'un \'el\'ement de $G$ sont soit
tous dans $\{ -1,1\}$, soit tous dans $\{ -i,i\}$.

La repr\'esentation naturelle de $G$ est irr\'eductible. Le centre de $G$
est d'ordre $4$. Le carr\'e d'un \'el\'ement de $G$ est une matrice
diagonale \`a coefficients diagonaux dans $\{ -1,1\}$. Ainsi l'exposant de
$G$ est $4$.

Etant donn\'ee une r\'eflexion g\'en\'eratrice $r$ de $\tW$,
on consid\`ere l'ensemble des r\'eflexions dans $G\cdot r$. C'est un ensemble
\`a $4$ \'el\'ements. Pour chacun de ses \'el\'ements, on consid\`ere
une forme lin\'eaire d\'efinissant l'hyperplan de r\'eflexion et
enfin, on note $p_r$ le produit de ces formes lin\'eaires.
On a

$$p_s=X_1X_2X_3X_4,\ \ 
p_t=((X_1+X_2)^2-(X_3+X_4)^2)((X_1-X_2)^2-(X_3-X_4)^2)$$
$$
p_u=(X_1^2+X_2^2)(X_3^2+X_4^2),\ \
p_v=(X_1^2-X_4^2)(X_2^2-X_3^2),\ \
p_w=(X_1^2-X_2^2)(X_3^2-X_4^2).
$$

L'alg\`ebre $\BC[X_1,X_2,X_3,X_4]^G$ est engendr\'ee par
$p_s,p_t,p_u,p_v,p_w$. C'est une hypersurface, \ie, le noyau
du morphisme canonique $\BC[Y_s,Y_t,Y_u,Y_v,Y_w]\to \BC[X_1,X_2,X_3,X_4]^G$,
$Y_r\mapsto p_r$, est engendr\'e par un polyn\^ome $R$.
L'expression de ce polyn\^ome est plus simple dans une nouvelle base
des invariants~:
$$
\begin{array}{cclcl}
p_1 &=&X_1^4+X_2^4+X_3^4+X_4^4   &=& -8p_s+p_t+3p_u+2p_v+p_w  \\
p_2 &=&2(X_1^2X_2^2 +X_3^2X_4^2) &=& p_u+2p_v+p_w  \\
p_3 &=&2(X_1^2X_3^2 +X_2^2X_4^2) &=& p_u+p_w  \\
p_4 &=&2(X_1^2X_4^2 +X_2^2X_3^2) &=& p_u-p_w  \\
p_5 &=&4X_1X_2X_3X_4             &=& 4p_s
\end{array}
$$

Alors, 
$$R=Y_1^2Y_5^2-2Y_1Y_2Y_3Y_4+Y_2^2Y_3^2+Y_2^2Y_4^2+Y_3^2Y_4^2
-Y_5^2(Y_2^2+Y_3^2+Y_4^2)+Y_5^4$$

L'action de $\tW$ sur $\BC[X_1,X_2,X_3,X_4]^G$ se factorise par
$W=\tW/G$. L'action induite de $W$ sur $\oplus_i \BC Y_i$
est la repr\'esentation de r\'eflexion de $W=\Sn_6$. Dans son action
sur $\BC^5$, le groupe $\Sn_6$ laisse invariante l'hypersurface $R=0$,
qui est une singularit\'e quotient, de groupe $G$.

\section{Construction de $\tW$ \`a partir de $W$ et $G$}~

\medskip

Dans les parties pr\'ec\'edentes, nous avons \'etudi\'e le cas d'un groupe de r\'eflexion $W$ 
obtenu comme quotient d'un groupe de r\'eflexion $\tW$ par un sous-groupe distingu\'e 
d'intersection compl\`ete $G$. Dans cette partie, nous nous int\'eressons au probl\`eme inverse, 
\`a savoir l'existence d'un groupe de r\'eflexion $\tW$ \'etant donn\'es $W$ et $G$. 
Sauf mention du contraire nous ne supposons plus que $k$ est de caract\'eristique nulle.

\bigskip

\subsection{Un rappel de g\'eom\'etrie alg\'ebrique} 
On appellera sch\'ema un sch\'ema noeth\'erien s\'epar\'e sur le corps $k$.
Rappelons l'existence de cl\^otures normales de rev\^etements~:

\begin{lemme}
\label{cloture}
Soit $Y$ un sch\'ema normal, $H$ un groupe fini agissant sur $Y$ et
$\pi:X\to Y$ un rev\^etement \'etale connexe. Soit $f:Y\to Y/H$
l'application quotient.

Alors, il existe un rev\^etement \'etale connexe
$p:Z\to X$ de $X$ et un groupe fini $G$ agissant sur $Z$ tels que
$f\pi p:Z\to Y/H$ est l'application quotient par $G$.
\end{lemme}

\begin{proof}
Soit $\overline{k(X)}$ une cl\^oture s\'eparable de $k(X)$.
Soit $K'$ la cl\^oture normale de l'extension $k(X)/k(Y)^H$ dans
$\overline{k(X)}$ et $Z$ la normalisation de $X$ dans $K'$. Pour montrer que
$Z$ satisfait aux propri\'et\'es du lemme, il faut montrer que le
rev\^etement $p:Z\to X$ est non ramifi\'e.

 Soit 
$K=\bigcup_{F} F$ o\`u $F$ d\'ecrit les extensions finies de 
$k(Y)$ contenues dans $\overline{k(X)}$ telles que la normalisation de
$Y$ dans $F$ n'est pas ramifi\'ee sur $Y$. L'extension
$K/k(Y)^W$ est normale.
Puisque $k(X)$ est un sous-corps de $K$, on d\'eduit
que $K'$ est un sous-corps de $K$, donc que $\pi p:Z\to Y$ est non
ramifi\'ee et finalement que $p$ est non ramifi\'ee.
\end{proof}

\begin{prop}
\label{relev}
Soit $X$ un sch\'ema lisse connexe simplement connexe sur $k$,
$G$ un groupe fini agissant sur $X$
et $W$ un groupe fini agissant sur $X/G$.
On suppose que le lieu de ramification de $G$ a codimension au moins $2$.

Alors, il existe un unique groupe fini $\tW$ agissant sur $X$, contenant
$G$ comme sous-groupe distingu\'e et tel que $\tW/G=W$.
\end{prop}

\begin{proof}
Le probl\`eme est de montrer que le rev\^etement $X\to (X/G)/W$ est
galoisien. Il suffit de le montrer pour un ouvert dense de $X$~:
nous allons le faire pour $U$ le compl\'ementaire du lieu de
ramification de $G$ dans $X$.

Notons que $U/G$ est lisse, puisque $U$ l'est et $U\to U/G$ est \'etale.
Puisque la codimension du compl\'ementaire de $U$ dans $X$ est au moins $2$ et que
$X$ est lisse et simplement connexe, la vari\'et\'e $U$ est simplement connexe
\cite[X, Corollaire 3.3]{SGA1}.
Notons $\pi$ l'application quotient $X\to X/G$.
En outre, le th\'eor\`eme de puret\'e de Zariski-Nagata 
\cite[X, Th\'eor\`eme 3.1]{SGA1} affirme que
le lieu de ramification de $\pi^{-1}((X/G)_{lisse})\to(X/G)_{lisse}$ est
vide, car s'il ne l'\'etait pas, il serait purement de codimension $1$.
Par cons\'equent, $U/G=(X/G)_{lisse}$ est stable par $W$.
Le lemme \ref{cloture}
permet de conclure que le rev\^etement $U\to (U/G)/W$ est galoisien.
\end{proof}

\bigskip

\subsection{Application aux groupes de r\'eflexion} 
Nous utiliserons de la proposition pr\'ec\'edente le corollaire suivant~:

\begin{cor}
\label{relevement}
Soit $\tV$ un espace vectoriel de dimension finie sur $k$ muni de sa graduation
standard, $G$ un sous-groupe fini de $GL(\tV)$
et $W$ un groupe fini agissant de mani\`ere gradu\'ee sur $\tV/G$.
On suppose que $G$ ne contient pas de r\'eflexion et que $k$ est
alg\'ebriquement clos de caract\'eristique nulle.

Alors, il existe un sous-groupe $\tW$ de $GL(\tV)$
contenant $G$ comme sous-groupe distingu\'e et tel que $\tW/G=W$.
\end{cor}

\begin{proof}
Puisqu'un espace affine sur un corps alg\'ebriquement clos de
caract\'eristique nulle est simplement connexe,
la proposition \ref{relev} assure l'existence du groupe $\tW$.

Soit $a\in \tV^*$ et $w\in\tW$. Puisque $\tW$ agit de mani\`ere gradu\'ee
sur $k[\tV]^G$, l'\'el\'ement $w(\prod_{g\in G}g(a))$ est homog\`ene de
degr\'e $|G|$.
Par cons\'equent, les \'el\'ements $wg(a)$ sont tous de degr\'e $1$.
Ainsi, $\tW$ agit de mani\`ere gradu\'ee sur $k[\tV]$ et pr\'eserve
le sous-espace $\tV^*$ de $k[\tV]$, donc $\tW$ est un sous-groupe de
$GL(\tV)$.
\end{proof}

\smallskip
\begin{rem}
\begin{itemize}
\item[(i)]
Si $G$ contient des r\'eflexions, le r\'esultat n'est plus vrai~: on le voit
en prenant $G=W$ groupe de r\'eflexion ($\tV/G$ \'etant identifi\'e \`a $\tV$).
\item[(ii)]
Le corollaire n\'ecessite $k$ alg\'ebriquement clos. La partie
\S \ref{Gordre2} fournit des cas o\`u $k=\BR$ et $\tW$ agit sur
$\tV\otimes_\BR \BC$ mais pas sur $\tV$.

\item[(iii)]
Le corollaire reste vrai pour $k$ alg\'ebriquement clos de
caract\'eristique $p$, \`a condition
que les ordres de $G$ et $W$ soient premiers \`a $p$, puisque le groupe
fondamental d'un espace affine sur un tel corps n'a pas de quotient non
trivial d'ordre premier \`a $p$~ (\'etant donn\'e trois groupes finis
$H\lhd H'\lhd \Gamma$ avec $H=O^{p'}(H)$ et $p\not|\ [\Gamma:H]$, alors
$H=O^{p'}(H')$, donc $H\lhd \Gamma$).

\item[(iv)]
Le corollaire est faux en g\'en\'eral pour $k$ alg\'ebriquement clos
de caract\'eristique $p$, comme le montre l'exemple suivant.
On prend $E$ espace vectoriel de dimension $2$ sur $k$ muni d'un endomorphisme
de Frobenius $F$.
Soit $B$ un sous-groupe de Borel $F$-stable de $GL(E)$
et $U$ le radical unipotent de $B$.
Soit $\tV=E/U^F$, $G=U^{F^2}/U^F$ et
$W=B^{F^2}/U^{F^2}$. Alors, $\tV$ est un espace affine
de dimension $2$ sur $k$, $\tV\to \tV/G$ est non ramifi\'ee, $W$ est
un groupe de r\'eflexion sur $\tV/G$. N\'eanmoins, le rev\^etement
$\tV\to (\tV/G)/W$ n'est pas galoisien ; une cl\^oture normale
est donn\'ee par $E\to \tV$.
\end{itemize}
\end{rem}

\begin{thm}
Soit $\tV$ un $k$-espace vectoriel muni de sa graduation standard,
$G$ un sous-groupe fini de $GL(\tV)$ et $Z=\tV/G$.
Soit $V$ un $k$-espace vectoriel gradu\'e,
$W$ un groupe de r\'eflexion sur $V$
et $L$ un sous-espace affine homog\`ene de $V/W$ tel qu'il existe un
isomorphisme $Z\iso L\times_{V/W}V$ compatible \`a l'action de $k^\times$.

On suppose que $G$ ne contient pas de r\'eflexion
et que $k$ est un corps alg\'ebriquement clos de caract\'eristique
nulle.

Alors, il existe un (unique) groupe de r\'eflexion $\tW$ sur
$\tV$ contenant $G$ comme sous-groupe distingu\'e et tel que $\tW/G=W$.
En outre, $G$ est un bon sous-groupe distingu\'e de $\tW$
et $W$ agit trivialement sur
le quotient de $V$ par l'espace tangent \`a $Z$ en $0$.
\end{thm}

\begin{proof}
Notons que $W$ agit fid\`element sur $Z$.
L'existence de $\tW$ r\'esulte alors du corollaire \ref{relevement}.
Puisque $\tV/\tW\simeq Z/W$ est un espace affine, $\tW$ est
de r\'eflexion.

Soit $E$ l'espace tangent \`a $Z$ en $0$, vu comme sous-espace de $V$.
Puisque $W$ agit fid\`element sur $E$ et qu'il est de r\'eflexion,
il agit trivialement sur $V/E$. Quitte \`a remplacer $V$ par $E$,
on est dans le cas o\`u $V$ est l'espace tangent \`a $Z$ en $0$ et on
conclut par le th\'eor\`eme \ref{ci} que $Z$ est une intersection
compl\`ete.
\end{proof}

\end{document}